\documentclass{article}
\usepackage[top=2.5cm, bottom=2.5cm, left=2.5cm, right=2.5cm]{geometry}
\usepackage{amsfonts,amsmath,amssymb,amsthm,amscd,mathrsfs,xspace,multicol}
\usepackage{latexsym}
\usepackage[all]{xy}
\usepackage[breaklinks=true]{hyperref}
\usepackage{graphicx}
\usepackage{epsfig}

\input{epsf}

\newcommand{\Set}{{\rm Set}}
\newcommand{\FinSet}{{\rm FinSet}}
\newcommand{\Vect}{{\rm Vect}}
\newcommand{\Cat}{{\rm Cat}}
\newcommand{\Grpd}{{\rm Grpd}}
\newcommand{\Hecke}{{\rm Hecke}}
\newcommand{\Span}{{\rm Span}}
\newcommand{\PermRep}{{\rm PermRep}}
\newcommand{\Nice}{{\rm Nice}}

\newcommand{\GSet}{G{\rm Set}}
\newcommand{\Cocont}{{\rm Cocont}}

\renewcommand{\hom}{{\rm hom}}
\newcommand{\To}{\Rightarrow}
\renewcommand{\to}{\rightarrow}
\newcommand{\maps}{\colon}

\newcommand{\SL}{\mathrm{SL}}
\newcommand{\C}{{\mathbb C}}
\newcommand{\F}{{\mathbb F}}

\newcommand{\Aut}{{\rm Aut}}
\newcommand{\Over}{/\!/}

\newtheorem{theorem}{Theorem}
\newtheorem*{theorem*}{Theorem}
\newtheorem{definition}[theorem]{Definition}
\newtheorem{lemma}[theorem]{Lemma}
\newtheorem{claim}[theorem]{Claim}
\newtheorem*{claim*}{Claim}

\newdir{ >}{{}*!/-7pt/@{>}}
\newdir{> >}{*!/0pt/@{>}*!/-5pt/@{>}}

\begin{document}
\sloppy
\title{The Hecke Bicategory}
\author{Alexander E.\ Hoffnung\\
       Department of Mathematics and Statistics\\ University of Ottawa\\
       Ottawa, ON, K1N 6S5 Canada\\
       email:\;hoffnung@uottawa.ca}
\maketitle

\begin{abstract}
We present an application of the program of groupoidification leading up to a sketch of a categorification of the Hecke algebroid --- the category of permutation representations of a finite group.  As an immediate consequence, we obtain a categorification of the Hecke algebra.  We suggest an explicit connection to new higher isomorphisms arising from incidence geometries, which are solutions of the Zamolodchikov tetrahedron equation.  This paper is expository in style and is meant as a companion to {\em Higher Dimensional Algebra VII: Groupoidification} and an exploration of structures arising in the work in progress, {\em Higher Dimensional Algebra VIII: The Hecke Bicategory}, which introduces the Hecke bicategory in detail.
\end{abstract}

\tableofcontents

\section{Introduction}
Categorification is, in part, the attempt to shed new light on familiar mathematical notions by replacing a set-theoretic interpretation with a category-theoretic analogue.  Loosely speaking, categorification replaces sets, or more generally $n$-categories, with categories, or more generally $(n+1)$-categories, and functions with functors.  By replacing {\em interesting} equations by isomorphisms, or more generally equivalences, this process often brings to light a new layer of structure previously hidden from view.  While categorification is not a systematic process --- in other words, finding this new layer of structure may require a certain amount of creativity --- the reverse process of {\em decategorification} should be a systematic way of recovering the original set-theoretic structure or concept.  The key idea is that considering a process of categorification requires, as a first step, a definition of the corresponding decategorification process.  We then think of categorification simply as a heuristic tool allowing us to `undo' the process of decategorification.

In {\em Higher Dimensional Algebra VII: Groupoidification}\;\cite{HDA7}, Baez, Walker and the author introduced a program called {\em groupoidification} initiated by Baez, Dolan, and Trimble, and aimed at categorifying various notions from linear algebra, representation theory and mathematical physics.  The very simple idea was that one could replace vector spaces by groupoids, i.e., categories with only isomorphisms, and replace linear operators by spans of groupoids.  In fact, what we really did was define a systematic process called {\em degroupoidification}:
\[\textrm{groupoids}\mapsto \textrm{vector spaces}\]
\[\textrm{spans of groupoids} \mapsto \textrm{matrices}\]
\noindent Thus, groupoidification is a form of categorification.  We then suggested some applications of groupoidification to Hall algebras, Hecke algebras, and Feynman diagrams, so that other researchers could begin to categorify {\em their} favorite notions from representation theory.

In this paper, we give an expository account of a theory of {\em categorified intertwining operators} or {\em categorified Hecke operators} for representations of a very basic type: the {\em permutation representations} of a finite group.  Following the description of categorification above, this suggests the study of a $2$-dimensional category-theoretic structure and a decategorification functor.  We describe this $2$-dimensional structure, which we call the {\em Hecke bicategory}, in Section\;\ref{two}.  Pairing the Hecke bicategory with the degroupoidification functor we are then able to state a categorification theorem as Claim\;\ref{main}.

The statements of the main results are as follows.  For each finite group $G$, there is an equivalence of categories, or more precisely algebroids, between the category of permutation representations of $G$ --- the {\em Hecke algebroid} of {\em Hecke operators} --- and the degroupoidification of the {\em Hecke bicategory} of {\em categorified Hecke operators}, which has finite $G$-sets as objects and is enriched over the monoidal bicategory of spans of groupoids.  In other words, {\em the Hecke bicategory categorifies the category of permutation representations}.  When $G$ is the simple Lie group over a finite field of $q$ elements attached to a Dynkin diagram, then one can choose a Borel subgroup $B$, and construct the $G$-set $X = G/B$, known as the flag complex.  The choice of one object $X$ in $\Hecke(G)$ yields a groupoid $\Hecke(X,X)$ and a span called composition.  {\em The groupoid $\Hecke(X,X)$ and accompanying span categorify the usual Hecke algebra for the chosen Dynkin diagram and prime power $q$}.

The term `Hecke algebra' is seen in several areas of mathematics.  The Hecke algebras we consider, the {\em Iwahori-Hecke algebras}, are one-parameter deformations of the group algebras of Coxeter groups.  In the theory of modular forms, or more generally, automorphic representations, Hecke algebras are commutative algebras of Hecke operators.  These operators can be expressed by means of double cosets in the modular group, or more generally, with respect to certain compact subgroups.  Here we use the term `Hecke operator' to highlight the relationship between the intertwining operators between permutation representations of a finite group and the Hecke operators acting on the modular group.  We discuss an example of Hecke operators in terms of ``flag-flag relations" in the setting of Coxeter groups in Section\;\ref{A2}.  To describe Hecke algebras, one may use relations between varieties of the form $G/H$ for various subgroups $H$, namely the {\em discrete} subgroups, or more generally, the {\em compact} subgroups.  So, we think of Hecke algebras as algebras of ``Hecke operators" in a slightly generalized sense.  We are changing the groups, and making them {\em finite}, so that instead of varieties $G/H$, we have certain finite sets $G/H$.  Thus we think of a Hecke algebra as the algebra of intertwining operators from a permutation representation to itself.  Generalizing this slightly, we think of intertwining operators between permutation representations in general as {\em Hecke operators}.

We make use of the techniques of groupoidification along with the machinery of enriched bicategories and some very basic topos theory.  Thus, this paper is intended to give an introduction to some concepts which should play a significant role as the subject of categorified representation theory continues to develop.  More detailed accounts of the necessary structures along with proofs will be presented in papers in progress by the author\;\cite{Hoffnung1},\cite{Hoffnung2}, Kenney and Pronk\;\cite{KePr} and by the author in collaboration with John Baez\;\cite{HDA8}.  We now proceed to give a brief overview and explanation of each section.

\subsection{Matrices and Spans}\label{Intro1}

In Section\;\ref{spans}, we give a heuristic discussion a very simple notion of categorification.  In particular, we recall the basic notions of {\em spans}, also known as {\em correspondences}, and see that {\em spans of finite sets categorify matrices of natural numbers}.  Decategorification can be defined using just set cardinality and the free vector space construction, and this process is indeed functorial, since composition of spans by pullback corresponds to matrix multiplication.

After discussing the example of linear operators, we pass to the intertwining operators for a finite group $G$.  Then we need to consider not spans of finite sets but spans of finite $G$-sets.  Since the maps in the spans are now $G$-equivariant, the corresponding matrices should also be $G$-equivariant.  This prompts us to recall the relationships between finite $G$-sets and permutation representations of $G$.  There is a faithful, essentially surjective functor from the category of $G$-sets to the category of permutation representations of $G$.  However, this functor is not {\em full}.  Spans in a category with pullbacks naturally form a bicategory.  Since $G$-sets and permutation representations of $G$ are closely related, except that there are `not enough' maps of $G$-sets, the bicategory of spans of $G$-sets is a first clue in constructing  categorified permutation representations.  We will return to the role of spans of $G$-sets in Sections\;\ref{spans}, \;\ref{alternate} and\;\ref{applications}.

\subsection{Groupoidification and Enriched Bicategories}

To understand groupoidification, we need to recall the construction and basic properties of the {\em degroupoidification functor} defined in \cite{HDA7}.  We discuss this functor and extend it to a functor on the bicategory of spans of groupoids in Section\;\ref{degroupoidification}.  Our intention is to define the Hecke bicategory as an enriched structure, but the replacement of functors in the monoidal $2$-category of groupoids by spans in the monoidal bicategory of spans of groupoids necessitates a generalization of enriched category theory to enriched bicategory theory.   Definitions of enriched bicategories were developed independently by the author and in the unpublished Ph.D. thesis of Carmody\;\cite{Car}, which is reproduced in part by Forcey\;\cite{For}.  We present a partial definition explaining the basic idea in this work and will define the full structure along with a {\em change of base} theorem in the forthcoming paper\;\cite{Hoffnung2}.  Change of base for enriched bicategories is analogous to the theorem for enriched categories, and together with the degroupoidification functor, allows us to obtain a categorification theorem for Hecke operators.

The degroupoidification functor takes the monoidal bicategory $\Span(\Grpd)$ of spans of finite groupoids to the monoidal category $\Vect$.  We briefly recall this functor here.  A groupoid is sent to the free vector space on its set of isomorphism classes of objects.  A span of groupoids is sent to a linear operator using the {\em weak} or {\em pseudo} pullback of groupoids and the notion of {\em groupoid cardinality}\;\cite{BaezDolan:2001}.  That is, we think of a span of groupoids as a categorified or groupoid-valued matrix in much the same way as we think of a span of sets as a set-valued matrix, where set cardinality is replaced by groupoid cardinality.

The notion of enriched bicategories is then used to make the description of our decategorification processes precise.  Given a monoidal bicategory $\mathcal{V}$, a $\mathcal{V}$-enriched bicategory consists of a set of objects, $\hom$-objects in $\mathcal{V}$, composition morphisms in $\mathcal{V}$, and further structure and axioms, all of which live in the monoidal bicategory $\mathcal{V}$.  The only theorem about enriched bicategories that we will need is a {\em change of base} theorem.  In particular, given a functor $\mathcal{F}\maps \mathcal{V} \to \mathcal{V'}$ and a $\mathcal{V}$-enriched bicategory, then we obtain a $\mathcal{V'}$-enriched bicategory with $\hom$-objects $f(\hom(x,y))$.  If $\mathcal{V}$ is a monoidal category, then a $\mathcal{V}$-enriched bicategory is a $\mathcal{V}$-enriched category in the traditional sense\;\cite{Kelly}.  Enriched bicategories and change of base are discussed in more detail in Section\;\ref{enriched}.

\subsection{The Hecke Bicategory}

The Hecke bicategory $\Hecke(G)$ of a finite group $G$ is a categorification of the permutation representations of $G$.  This is the main construction of this paper.  Section\;\ref{two} gives a the basic structure of this family of enriched bicategories.  For each finite group $G$, $\Hecke(G)$ is an enriched bicategory over the monoidal bicategory of spans of groupoids.  The structure of this monoidal bicategory is given roughly in Section\;\ref{degroupoidification}.
A more complete and general account of monoidal bicategories (and monoidal tricategories) of spans will be given in\;\cite{Hoffnung1}.

The Hecke bicategory is constructed to study the Hecke operators between permutation representations $X$ and $Y$ of a finite group $G$ by keeping track of information about the $G$-orbits in $X\times Y$ as the $\hom$-groupoid $(X\times Y)\Over G$.  This double-slash notation denotes the {\em action groupoid}, which we recall in Definition\;\ref{actiondefn}.  The composition process between such groupoids is closely related to the {\em push-tensor-pull} construction familiar from geometric representation theory.

Now, applying the change of base theorem of enriched bicategories together with the degroupoidification functor, we obtain a $\Vect$-enriched category from $\Hecke(G)$.  The resulting $\Vect$-enriched category is equivalent to the $\Vect$-enriched category of permutation representations of $G$.  This is our main theorem and is stated in Claim\;\ref{main}.

\subsection{Spans of Groupoids and Cocontinuous Functors}

In Section\;\ref{algebroid}, the bicategory of spans of $G$-sets appear in our study of the category of permutation representations.  We view spans of $G$-sets as categorified $G$-equivariant matrices, but do not specify a decategorification process in this setting, although we discuss the importance of set cardinality.  This bicategory of spans $\Span(\GSet)$ plays two important roles in our attempt to understand categorified representation theory, which we briefly discuss here.

In Section\;\ref{A2}, we will explain the categorification of the Hecke algebra associated to the $A_2$ Dynkin diagram as a $\hom$-category in $\Span(\GSet)$.  This makes the categorification explicit in that we see the Hecke algebra relations holding {\em up to isomorphism} with special spans of $G$-sets acting as {\em categorified generators}.  The isomorphisms representing the relations come explicitly from incidence geometries --- in this example, projective plane geometry --- associated to the Dynkin diagram.

To facilitate this point of view, we describe a monoidal functor in Section\;\ref{monoidal} that takes a groupoid to the corresponding presheaf category and linearizes a span of groupoids to a cocontinuous functor by a {\em push-tensor-pull} process familiar from many geometric constructions in representation theory, but quite generally applicable in the theory of Grothendieck toposes.  This is a functor $\mathcal{L}\maps \Span\to \Cocont$, where $\Cocont$ is the monoidal $2$-category of presheaf categories on groupoids, cocontinuous functors, and natural isomorphisms.  By change of base, we obtain a bicategory enriched over certain presheaf toposes, the objects of which are interpreted as spans of $G$-sets in Section\;\ref{presheaf}.

We denote the functor taking a groupoid to its category of presheaves $\mathcal{L}$ suggesting that this might be interpreted as {\em categorified linearization}. In this interpretation one thinks of a groupoid as a {\em categorified vector space} that is equipped with a chosen basis --- its set of isomorphism classes of objects or connected components.  From a Grothendieck topos that is equivalent to a category of presheaves on a groupoid, we can always recover the groupoid.  Where a basis can always be recovered up to isomorphism from a vector space, a groupoid can be recovered up to equivalence from such a presheaf category.  Then we think of such a topos as an {\em abstract categorified vector space}.

The relationship between categorified representation theory enriched over monoidal bicategories of spans of spaces and enrichment over corresponding $2$-categories of sheaves on the spaces and cocontinuous functors between these will be the focus of future work.  From our point of view, these monoidal $2$-categories are the basic objects of study in categorified linear algebra, and in this sense, {\em categorified representation theory is enriched over categorified linear algebra}.  Part of the motivation for this line of work is to better understand and unify homology theories which arise in geometric representation theory.

\subsection{The Categorified Hecke Algebra and Zamolodchikov Equation}

As already mentioned, Section\;\ref{applications} discusses the main corollary, Claim\;\ref{maincorollary} --- a categorification of the Hecke algebra of a Coxeter group --- as well as possible future directions in low-dimensional topology and higher-category theory.  In Section\;\ref{hecke}, we recall the notion of Hecke algebras associated to Dynkin diagrams and prime powers.  We describe how a categorification of the Hecke algebra naturally arises from the Hecke bicategory.  

Finally, in Section\;\ref{A2}, we describe a concrete example of the categorified Hecke algebra in terms of spans of $G$-sets.  We describe solutions to the Zamolodchikov tetrahedron equation, which we hope will lead to constructions of braided monoidal $2$-categories as pointed out by Kapranov and Voevodsky\;\cite{KV}, and eventually tangle invariants.

Other approaches to categorified Hecke algebras and their
representations have been studied by a number of authors, building on
Kazhdan--Lusztig theory \cite{KaLu}.  One key step was Soergel's
introduction of what are nowadays called Soergel bimodules
\cite{Rou,Soe}.  Also important was Khovanov's categorification
of the Jones polynomial \cite{Kh} and work by Bernstein, Frenkel,
Khovanov and Stroppel on categorifying Temperley--Lieb algebras,
which are quotients of the type $A$ Hecke algebras
\cite{BFK,Str2}.  A diagrammatic interpretation of the Soergel
bimodule category was developed by Elias and Khovanov \cite{ElKh}, and
a geometric approach led Webster and Williamson \cite{WW} to deep
applications in knot homology theory.  This geometric interpretation
can be seen as going beyond the simple form of groupoidification we
consider here, and considering groupoids in the category of schemes.

\section{Matrices, Spans and $G$-Sets}\label{spans}
\subsection{Spans as Matrices}\label{matrix}
The first tool of representation theory is linear algebra.  Vector spaces and linear operators have nice properties, which allow representation theorists to extract a great deal of information about algebraic gadgets ranging from finite groups to Lie groups to Lie algebras and their various relatives and generalizations.  We start at the beginning, considering the representation theory of finite groups.
Noting the utility of linear algebra in representation theory, this paper is fundamentally based on the idea that the heavy dependence of linear algebra on fields, very often the complex numbers, may at times obscure the combinatorial skeleton of the subject.  Then, we hope that by peeling back the soft tissue of the continuum, we may expose and examine the bones, revealing new truths by working directly with the underlying combinatorics.  In this section, we consider the notion of spans of sets, a very simple idea, which is at the heart of categorified representation theory.

A {\bf span of sets} from $X$ to $Y$ is a pair of functions with a common domain like so:
\[
\xymatrix{
& S\ar[dr]^p\ar[dl]_q&\\
Y&&X
}
\]
\noindent We will often denote a span by its apex, when no confusion is likely to arise.

A span of sets can be viewed as a a matrix of sets
\[
\begin{picture}(180,150)(0,25)
 \includegraphics[scale=0.5]{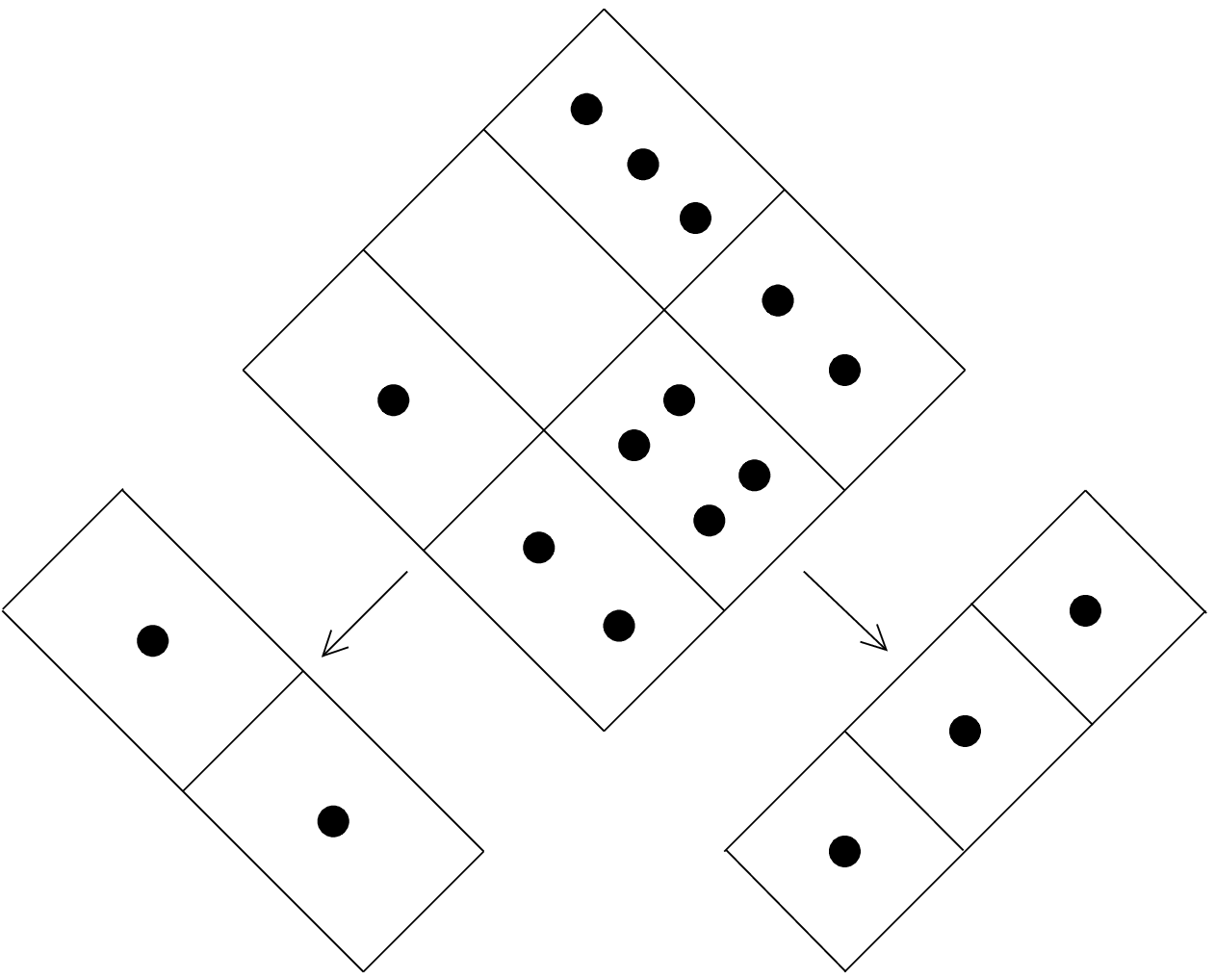}
 \put(-138,88){$q$}
 \put(-52,88){$p$}
 \put(-8,98){$X$}
 \put(-180,98){$Y$}
 \put(-80,170){$S$}
\end{picture}
\]
\noindent For each pair $(x,y)$, we have a set $S_{x,y} = p^{-1}(x)\cap q^{-1}(y)$.  In particular, if all the sets $S_{x,y}$ are finite, we can apply a process of `decategorification' to obtain a matrix of natural numbers $|S_{x,y}|$ --- a very familiar object in linear algebra.  In this sense, a span is a `categorification' of a matrix.  We will consider only {\em finite} sets throughout this paper.

In addition to spans giving rise to matrices, composition of spans gives rise to matrix multiplication.  Given a pair of composable spans
\[
\xymatrix{
& T\ar[dl]\ar[dr]^{q} && S\ar[dl]_{p}\ar[dr] & \\
Z && Y && X
}
\]
\noindent we define the composite to be the {\bf pullback} of the pair of functions $p\maps S\to Y$ and $q\maps T\to Y$, which is a new span
\[
\xymatrix{
&& TS\ar[dl]\ar[dr] &&\\
& T\ar[dl]\ar[dr]^{q} && S\ar[dl]_{p}\ar[dr] & \\
Z && Y && X
}
\]

\noindent where $TS$ is the subset of $T\times S$
\[\{(t,s)\subseteq T\times S\mid p(s)=q(t)\}\]
\noindent with the obvious projections to $S$ and $T$.  It is straightforward to check that this process agrees with matrix multiplication after decategorification.

In the above example, we turn spans of finite sets into matrices of natural numbers simply by counting the number of elements in each set $S_{x,y}$.  We note that there is a standard basis of the vector space of linear maps, and each basis element can be realized as a span of finite sets.  Thus, we can recover the entire vector space of linear maps by constructing the free vector space on this basis.  Further, we can turn a set $X$ into a vector space by constructing the free vector space with basis $X$.  Checking that composition of spans and matrix multiplication agree after taking the cardinality of the set-valued entries is the main step in showing that our decategorification process --- from the bicategory of spans of sets to the category of linear operators --- is {\em functorial}.

\subsection{Permutation Representations}\label{permrep}

Again we start with a very simple idea.  We want to study the actions of a finite group $G$ on finite sets --- {\em finite $G$-sets}.  These extend to {\em permutation representations of $G$}.  We fix the field of complex numbers and consider only complex vector spaces throughout this paper.

\begin{definition}
A {\bf permutation representation} of a finite group $G$ is a finite-dimensional representation of $G$ together with a chosen basis such that the action of $G$ maps basis vectors to basis vectors.
\end{definition}

\begin{definition}
An {\bf intertwining operator} $f\maps V\to W$ between permutation representations of a finite group $G$ is a linear operator from $V$ to $W$ that is $G$-equivariant, i.e., commutes with the actions of $G$.
\end{definition}

\noindent Finite $G$-sets can be linearized to obtain permutation representations of $G$.  In fact, this describes a relationship between the objects of the {\em category} of finite $G$-sets and the objects of the {\em category} of permutation representations of $G$. Given a finite group $G$, the category of finite $G$-sets has
\begin{itemize}
\item finite $G$-sets as objects,
\item $G$-equivariant functions as morphisms,
\end{itemize}
\noindent and the category of permutation representations $\PermRep(G)$ has
\begin{itemize}
\item permutation representations of $G$ as objects,
\item intertwining operators as morphisms.
\end{itemize}

\noindent One usually wants the morphisms in a category to preserve the structure on the objects.  Of course, an intertwining operator does not necessarily preserve the chosen basis of a permutation representation.  We can reconcile our choice of intertwining operators as morphisms, by noticing that $\PermRep(G)$ is equivalent to the category where each object is a finite-dimensional representation of $G$ with the {\em property} that there {\em exists} a basis preserved by the action of $G$.  So, we are justified in working with this version of $\PermRep(G)$.

A primary goal of this paper is to categorify the $q$-deformed versions of the group algebras of Coxeter groups known as Hecke algebras.  Of course, an algebra is a $\Vect$-enriched category with exactly one object, and the Hecke algebras are isomorphic to certain one-object subcategories of the $\Vect$-enriched category of permutation representations.  Thus, we refer to the category $\PermRep(G)$ as the {\em Hecke algebroid} --- a many-object generalization of the Hecke algebra.  We will construct a bicategory --- or more precisely, an {\em enriched bicategory} --- called the {\em Hecke bicategory} that categorifies the Hecke algebroid for any finite group $G$.

There is a functor from finite $G$-sets to permutation representations of $G$.  A finite $G$-set $X$ is linearized to a permutation representation $\widetilde{X}$, which is the free vector space on $X$.  As stated above, the maps between $G$-sets are $G$-equivariant functions --- that is, functions between $G$-sets $X$ and $Y$ that respect the actions of $G$.  Such a function $f\maps X\to Y$ gives rise to a $G$-equivariant linear map (or intertwining operator) $\widetilde{f} \maps \widetilde{X} \to \widetilde{Y}$.  However, there are many more intertwining operators from $\widetilde{X}$ to $\widetilde{Y}$ than there are $G$-equivariant maps from $X$ to $Y$.  In particular, the former is a complex vector space, while the latter is a finite set.  For example, an intertwining operator $\widetilde{f} \maps \widetilde{X} \to \widetilde{Y}$ may take a basis vector $x\in \widetilde{X}$ to any $\C$-linear combination of basis vectors in $\widetilde{Y}$, whereas a map of $G$-sets does not have the freedom of scaling or adding basis elements.

So, in the language of category theory the process of linearizing finite $G$-sets to obtain permutation representations is a faithful, essentially surjective functor, which is not at all full.

\subsection{Spans of G-Sets}\label{algebroid}

In the previous section, we discussed the relationship between finite $G$-sets and permutation representations.  In Section\;\ref{matrix}, we saw the close relationship between spans of finite sets and matrices of natural numbers.  There is an analogous relationship between spans of finite $G$-sets and $G$-equivariant matrices of natural numbers.

A {\bf span of finite $G$-sets} from a finite $G$-set $X$ to a finite $G$-set $Y$ is a pair of maps with a common domain
\[
\xymatrix{
& S\ar[dr]^p\ar[dl]_q&\\
Y&&X
}
\]
\noindent where $S$ is a finite $G$-set, and $p$ and $q$ are $G$-equivariant maps with respective codomains $X$ and $Y$.

Spans of finite $G$-sets are natural structures to consider in categorified representation theory because spans of sets --- and, similarly $G$-sets --- naturally form a bicategory.

The development of bicategories by Benabou\;\cite{Benabou} is an early example of categorification.  A (small) category consists of a {\em set of objects} and a {\em set of morphisms}.  A bicategory is a categorification of this concept, so there is a new layer of structure\;\cite{Leinster1}.  In particular, a (small) bicategory $\mathcal{B}$ consists of:

\begin{itemize}
\item a set of objects $x,y,z\ldots$,
\item for each pair of objects, a set of morphisms,
\item for each pair of morphisms, a set of $2$-morphisms,
\end{itemize}
\noindent and given any pair of objects $x,y$, this data forms a $\hom$-category $\hom(x,y)$ which has:
\begin{itemize}
\item $1$-morphisms $x \to y$ of $\mathcal{B}$ as objects,
\item $2$-morphisms:
\[
\xy
(-11,0)*+{x}="2";
(11,0)*+{y}="3";
(-8,0)*+{\bullet}="4";
(8,0)*+{\bullet}="6";
{\ar@/^1.65pc/ "4";"6"};
{\ar@/_1.65pc/ "4";"6"};
{\ar@{=>} (0,3)*{};(0,-3)*{}} ;
\endxy
\]
\noindent of $\mathcal{B}$ as morphisms.
\end{itemize}
There is a {\em vertical composition} of $2$-morphisms,
\[
\xy
(-11,0)*+{}="2";
(11,0)*+{}="3";
(-8,0)*+{\bullet}="4";
(8,0)*+{\bullet}="6";
{\ar@/^1.65pc/ "4";"6"};
{\ar "4";"6"};
{\ar@/_1.65pc/ "4";"6"};
{\ar@{=>} (0,5)*{};(0,1)*{}};
{\ar@{=>} (0,-1)*{};(0,-5)*{}};
\endxy
\]
\noindent as well as a {\em horizontal composition},
\[
\xy
(-22,0)*+{}="2";
(0,-3)*+{}="3";
(22,0)*+{}="5";
(-19,0)*+{\bullet}="4";
(0,0)*+{\bullet}="6";
(19,0)*+{\bullet}="7";
{\ar@/^1.65pc/ "4";"6"};
{\ar@/_1.65pc/ "4";"6"};
{\ar@/^1.65pc/ "6";"7"};
{\ar@/_1.65pc/ "6";"7"};
{\ar@{=>} (-10,3)*{};(-10,-3)*{}};
{\ar@{=>} (10,3)*{};(10,-3)*{}};
\endxy
\]
\noindent and these are required to satisfy certain coherence axioms, which make these operations simultaneously well-defined.

Benabou's definition followed from several important examples of bicategories, which he presented in\;\cite{Benabou}, and which are very familiar in categorified and geometric representation theory.  The first example is the bicategory of spans of sets, which has:
\begin{itemize}
\item sets as objects,
\item spans of sets as morphisms,
\item maps of spans of sets as $2$-morphisms.
\end{itemize}
\noindent We defined spans of sets in Section\;\ref{matrix}.  A {\bf map of spans of sets} from a span $S$ to a span $T$ is a function $f\maps S\to T$ such that the following diagram commutes:
\[
\xymatrix{
& S\ar[dr]^p\ar[dl]_q\ar[dd]^f&\\
Y&&X\\
& S'\ar[ul]^{q'}\ar[ur]_{p'} &
}
\]

For each finite group $G$, there is a closely related bicategory $\Span(\GSet)$, which has:
\begin{itemize}
\item finite $G$-sets as objects,
\item spans of finite $G$-sets as morphisms,
\item maps of spans of finite $G$-sets as $2$-morphisms.
\end{itemize}
\noindent The definitions are the same as in the bicategory of spans of sets, except for the obvious finiteness condition and that every arrow should be $G$-equivariant.

While this bicategory is a good candidate for a categorification of the Hecke algebroid $\PermRep(G)$, the theory of groupoidification allows for a related, but from our point of view, heuristically nicer categorification.  In what follows, we develop the necessary machinery to present this categorification, the Hecke bicategory $\Hecke(G)$, in the context of groupoidification.  In Sections\;\ref{alternate} and \;\ref{A2}, we return to $\Span(\GSet)$ and in future work we will make precise the relationship between spans of $G$-sets and the Hecke bicategory $\Hecke(G)$.

\section{Groupoidification and Enriched Bicategories}\label{groupoidification}

The following sections introduce the necessary machinery to present the Hecke bicategory, a categorification of the Hecke algebroid.  Enriched bicategories are developed for use in Section\;\ref{fundamentalsection} to construct the Hecke bicategory and state the main result, Theorem\;\ref{main}, and in Section\;\ref{alternate} to make a connection with the bicategory $\Span(\GSet)$ of Section\;\ref{algebroid}.

\subsection{Action Groupoids and Groupoid Cardinality}\label{action}

In this section, we draw a connection between $G$-sets and groupoids via the `action groupoid' construction.  We then introduce {\em groupoid cardinality}, which is the first step in describing the degroupoidification functor in the next section.

For any $G$-set, there exists a corresponding groupoid, called the {\em action groupoid} or {\em weak quotient}:

\begin{definition}\label{actiondefn}
Given a group $G$ and a $G$-set $X$, the {\bf action groupoid} $X\Over G$ is the category which has:
\begin{itemize}
\item elements of $X$ as objects,
\item pairs $(g,x)\in G\times X$ as morphisms $(g,x)\maps x \to x'$, where $g\cdot x = x'$.
\end{itemize}
\noindent Composition of morphisms is defined by the product in $G$.
\end{definition}

\noindent Of course, associativity follows from associativity in $G$ and the construction defines a groupoid since any morphism $(g,x)\maps x\to x'$ has an inverse $(g^{-1},x')\maps x' \to x$.

So every finite $G$-set defines a groupoid, and we will see in Section\;\ref{two} that the weak quotient of $G$-sets plays an important role in understanding categorified permutation representations.

Next, we recall the definition of groupoid cardinality\;\cite{BaezDolan:2001}:

\begin{definition}
Given a (small) groupoid $\mathcal{G}$, its {\bf groupoid cardinality} is defined as:
\begin{equation*}
\label{eq:groupoid_cardinality}
      |\mathcal{G}| = \sum_{\rm isomorphism\; classes \; of \; objects\; [x]}
\frac{1}{|\Aut(x)|}
\end{equation*}
\noindent If this sum diverges, we say $|\mathcal{G}| = \infty$.
\end{definition}

\noindent In this paper, we will only consider {\em finite groupoids} --- groupoids with a finite set of objects and finite set of morphisms. In general, we could allow groupoids with infinitely many isomorphism classes of objects, and the cardinality of a groupoid would take a value in the non-negative real numbers in case the sum converges.  Generalized cardinalities have been studied by a number of authors\;\cite{FL,Kim,Leinster2,Weinstein}.

We can think of groupoid cardinality as a form of categorified division analogous to the quotient of a $G$-set by its action of $G$ in the case when this action is free.  See the paper of Baez and Dolan\;\cite{BaezDolan:2001}.  In particular, we have the following equation:
\[ |X\Over G| = |X|/|G| \]
whenever $G$ is a finite group acting on a finite set $X$.

In the next section, we define the degroupoidification using the notion of groupoid cardinality.

\subsection{Degroupoidification}\label{degroupoidification}

In this section, we recall some of the main ideas of groupoidification.  Of course, in practice this means we will discuss the corresponding process of decategorification --- the degroupoidification functor.

To define degroupoidification in\;\cite{HDA7}, a functor was constructed from the category of spans of groupoids to the category of linear operators between vector spaces.  In the present setting, we will need to extend degroupoidification to a functor between bicategories.

We extend the functor to a bicategory $\Span(\Grpd)$ which has:
\begin{itemize}
\item (finite) groupoids as objects;
\item spans of (finite) groupoids as $1$-morphisms;
\item `isomorphism classes of equivalences' of spans of (finite) groupoids as $2$-morphisms.
\end{itemize}

Since all groupoids that show up in this paper arise from the action groupoid construction on finite $G$-sets, there is no problem restricting our attention to finite groupoids.

Arbitrary spans of groupoids should form a tricategory, which has not only {\em maps of spans} as $2$-morphisms, but also {\em maps of maps of spans} as $3$-morphisms.  For our purposes we restrict this structure to a bicategory.  While there are more sophisticated ways of obtaining such a bicategory, we do so by taking {\em isomorphism classes of equivalences of spans} as $2$-morphisms\;\cite{Hoffnung1}.

Spans of finite groupoids are categorified matrices of non-negative rational numbers in the same way that spans of finite sets are categorified matrices of natural numbers.  A {\em span of groupoids} is a pair of functors with common domain, and we can picture one of these roughly as follows:

\[
\begin{picture}(230,145)(0,30)
 \includegraphics[scale=0.5]{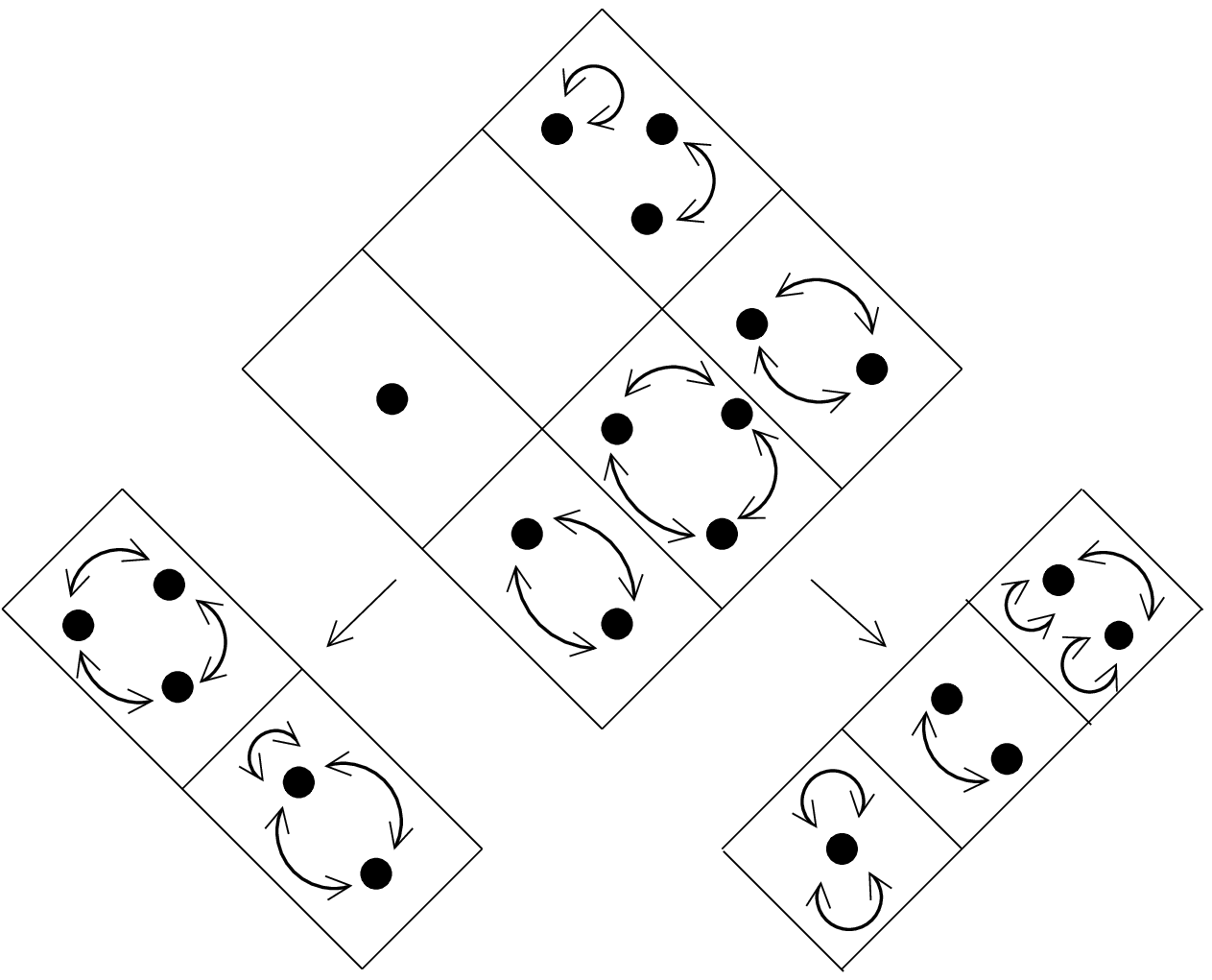}
 \put(-138,88){$q$}
 \put(-52,88){$p$}
 \put(-8,98){$\mathcal{G}$}
 \put(-180,98){$\mathcal{H}$}
 \put(-80,170){$\mathcal{S}$}
\end{picture}
\]

\noindent Whereas one uses set cardinality to realize spans of sets as matrices, we can use groupoid cardinality to obtain a matrix from a span of groupoids.

We have seen evidence that a span of groupoids is a categorified matrix, so we will want to think of a groupoid as a categorified vector space.  To make these notions precise, we recall the degroupoidification functor
\[\mathcal{D}\maps\Span\to \Vect.\]

Given a groupoid $\mathcal{G}$, we obtain a vector space $\mathcal{D}(\mathcal{G})$, called the {\bf degroupoidification} of $\mathcal{G}$, by taking the free vector space on the set of isomorphism classes of objects of $\mathcal{G}$.

We say a groupoid $\mathcal{V}$ {\em over} a groupoid $\mathcal{G}$
\[
\xymatrix{
\mathcal{V}\ar[d]^{p}\\
\mathcal{G}
}
\]
\noindent is a {\bf groupoidified vector}.  In particular, from the functor $p$ we can produce a vector in $\mathcal{D}(\mathcal{G})$ in the following way.

The {\bf full inverse image} of an object $x$ in $\mathcal{G}$ is the groupoid $p^{-1}(x)$, which has:
\begin{itemize}
\item objects $v$ of $\mathcal{V}$, such that $p(v)\cong x$, as objects,
\item morphisms $v\to v'$ in $\mathcal{V}$ as morphisms.
\end{itemize}
\noindent We note that this construction depends only on the isomorphism class of $x$.  Since the set of isomorphism classes of $\mathcal{G}$ determine a basis of the corresponding vector space, the vector determined by $p$ can be defined as
\begin{equation*}
      \sum_{\rm isomorphism\; classes \; of \; objects\; [x]}
|p^{-1}(x)|[x],
\end{equation*}
\noindent where $|p^{-1}(x)|$ is the groupoid cardinality of $p^{-1}(x)$.  We note that a `groupoidified basis' can be obtained in this way as a set of functors from the terminal groupoid $\mathbf{1}$ to representative objects of each isomorphism class of $\mathcal{G}$.  A {\bf groupoidified basis} of $\mathcal{G}$ is a set of groupoids $\mathcal{V} \to \mathcal{G}$ over $\mathcal{G}$ such that the corresponding vectors give a basis of the vector space $\mathcal{D}(\mathcal{G})$.

Given a span of groupoids
\[
\xymatrix{
& \mathcal{S}\ar[dl]_{q}\ar[dr]^{p} &\\
\mathcal{H}&&\mathcal{G}
}
\]
\noindent we want to produce a linear map $\mathcal{D}(\mathcal{S}) \maps \mathcal{D}(\mathcal{G})\to\mathcal{D}(\mathcal{H})$.  The details are checked in\;\cite{HDA7}.  Here we show only that given a basis vector of $\mathcal{D}(\mathcal{G})$ viewed as a groupoidified basis vector of $\mathcal{G}$, the span $\mathcal{S}$ determines a vector in $\mathcal{D}(\mathcal{H})$.  To do this, we need the notion of the weak pullback of groupoids --- a categorified version of the pullback of sets.

Given a diagram of groupoids
\[
\xymatrix{
\mathcal{H}\ar[dr]_{q} & & \mathcal{G}\ar[dl]^{p}\\
&\mathcal{I}&
}
\]
\noindent the {\bf weak pullback} of $p\maps\mathcal{G}\to\mathcal{I}$ and $q\maps \mathcal{H}\to\mathcal{I}$ is the diagram
\[
\xymatrix{
&\mathcal{HG}\ar[dr]\ar[dl]&\\
\mathcal{H}\ar[dr]_{q} & & \mathcal{G}\ar[dl]^{p}\\
&\mathcal{I}&
}
\]
\noindent where $\mathcal{HG}$ is a groupoid whose objects are triples $(h,g,\alpha)$ consisting of an object $h \in \mathcal{H}$, an object $g \in \mathcal{G}$, and an isomorphism $\alpha \maps p(g) \to q(h)$ in $\mathcal{I}$.  A morphism in $\mathcal{HG}$ from $(h,g,\alpha)$ to $(h',g',\alpha')$ consists of a morphism
$f \maps g \to g'$ in $\mathcal{G}$ and a morphism $f' \maps h \to h'$ in $\mathcal{H}$ such that the following square commutes:
\[
\xymatrix{
p(g) \ar[d]_{p(f)} \ar[r]^{\alpha} & q(h) \ar[d]^{q(f')} \\
p(g') \ar[r]_{\alpha'} & q(h')
}
\]
\noindent As in the case of the pullback of sets, the maps out of $\mathcal{HG}$ are the obvious projections.  Further, this construction satisfies a certain universal property.  See\;\cite{Weber}, for example.

Now, given our span and a chosen {\em groupoidified basis vector}:
\[
\xymatrix{
&\mathcal{S}\ar[dr]\ar[dl] && \mathbf{1}\ar[dl]\\
\mathcal{H} &&\mathcal{G} &
}
\]
\noindent we obtain a groupoid over $\mathcal{H}$ by constructing the weak pullback:
\[
\xymatrix{
&&\mathcal{S}\mathbf{1}\ar[dr]\ar[dl] &\\
&\mathcal{S}\ar[dr]\ar[dl] && \mathbf{1}\ar[dl]\\
\mathcal{H} &&\mathcal{G} &
}
\]
\noindent Now, $\mathcal{S}\mathbf{1}$ is a groupoid over $\mathcal{H}$, and we can compute the resulting vector.  In general, to guarantee that the process of passing a groupoidified vector across a span defines a linear operator, we need to restrict to the so-called `tame' spans defined in\;\cite{HDA7}.  However, spans of finite groupoids are automatically tame, so we can safely ignore this issue.

One can check that the process described above defines a linear operator from a span of groupoids, and, further, that this process is functorial\;\cite{HDA7}. This is the degroupoidification functor.  Since isomorphic spans are sent to the same linear operator, it is straightforward to extend this to our bicategory of spans of groupoids by adding identity $2$-morphisms to the category of vector spaces and sending all $2$-morphisms between spans of groupoids to the corresponding identity $2$-morphism.

In the next section, we give the basics of the notion of {\em enriched bicategories}. We will see that constructing an enriched bicategory requires having a monoidal bicategory in hand.  The bicategory $\Span(\Grpd)$ defined above is, in fact, a monoidal bicategory --- that is, $\Span(\Grpd)$ has a tensor product, which is a functor
\[\otimes \maps \Span(\Grpd)\times \Span(\Grpd) \to \Span(\Grpd),\]
\noindent along with further structure and satisfying some coherence relations.  The structure of $\Span(\Grpd)$, or more generally, monoidal bicategories of spans in bicategories, will be described in more detail in\;\cite{Hoffnung1}.

We describe the main components of the tensor product on $\Span(\Grpd)$.  Given a pair of groupoids $\mathcal{G},\mathcal{H}$, the tensor product $\mathcal{G}\times\mathcal{H}$ is the product in $\Grpd$.  Further, for each pair of pairs of groupoids $(\mathcal{G},\mathcal{H}),(\mathcal{I},\mathcal{J})$ there is a functor:
\[ \otimes \maps \Span(\Grpd)(\mathcal{G}, \mathcal{H})\times \Span(\Grpd)(\mathcal{I},\mathcal{J}) \to \Span(\Grpd)(\mathcal{G}\times \mathcal{I}, \mathcal{H}\times \mathcal{J}) ,\]
\noindent defined roughly as follows:
\[
\xy
  	    (-45,20)*+{\mathcal{S}}="1";
  	    (-55,10)*+{\mathcal{H}}="2";
  	    (-35,10)*+{\mathcal{G}}="3";
  	    (-15,20)*+{\mathcal{T}}="4";
  	    (-25,10)*+{\mathcal{J}}="5";
  	    (-5,10)*+{\mathcal{I}}="6";   	     	      	
  	    (30,20)*+{\mathcal{S}\times \mathcal{T}}="7";
  	    (20,10)*+{\mathcal{H}\times \mathcal{J}}="8";
  	    (40,10)*+{\mathcal{G}\times \mathcal{I}}="9";
  	    (-1,15)*+{}="10";
  	    (9,15)*+{}="11";
  	    (-30,15)*+{,}="12";
  	    (-1,-10)*+{}="14";
  	    (9,-10)*+{}="15";
  	    (-30,15)*+{,}="16";
  	    (-30,-10)*+{,}="27";  	
  	    (-45,0)*+{\mathcal{S}}="17";
  	    (-45,-20)*+{\mathcal{S}'}="13";  	
  	    (-55,-10)*+{\mathcal{H}}="18";
  	    (-35,-10)*+{\mathcal{G}}="19";
   	    (-15,0)*+{\mathcal{T}}="26";
  	    (-15,-20)*+{\mathcal{T}'}="20";
  	    (-25,-10)*+{\mathcal{J}}="21";
  	    (-5,-10)*+{\mathcal{I}}="22";   	     	      	
  	    (30,0)*+{\mathcal{S}\times \mathcal{T}}="23";
  	    (30,-20)*+{\mathcal{S}'\times \mathcal{T}'}="27";  	
  	    (20,-10)*+{\mathcal{H}\times \mathcal{J}}="24";
  	    (40,-10)*+{\mathcal{G}\times \mathcal{I}}="25";			
	      {\ar@{->}_{v} "4";"5"};
	      {\ar@{->}^{u} "4";"6"};
	      {\ar@{->}_{q} "1";"2"};
	      {\ar@{->}^{p} "1";"3"};
	      {\ar@{->}_{q\times v} "7";"8"};
	      {\ar@{->}^{p\times u} "7";"9"};
	      {\ar@{|->}_{} "10";"11"};
	      {\ar@{->}_{v} "26";"21"};
	      {\ar@{->}^{u} "26";"22"};
	      {\ar@{->}_{q} "17";"18"};
	      {\ar@{->}^{p} "17";"19"};
	      {\ar@{->}_{q\times v} "23";"24"};
	      {\ar@{->}^{p\times u} "23";"25"};
	      {\ar@{->}^{v'} "20";"21"};
	      {\ar@{->}_{u'} "20";"22"};
	      {\ar@{->}^{q'} "13";"18"};
	      {\ar@{->}_{p'} "13";"19"};
	      {\ar@{->} "17";"13"};
	      {\ar@{->} "26";"20"};
	      {\ar@{->}^{q'\times v'} "27";"24"};
	      {\ar@{->}_{p'\times u'} "27";"25"};
          {\ar@{|->}_{} "14";"15"};	
	      {\ar@{->}_{} "23";"27"};	
        {\ar@{=>}^<<{\scriptstyle \nu} (-51,-8); (-49,-12)};
        {\ar@{=>}_<<{\scriptstyle \mu} (-39,-8); (-41,-12)};
        {\ar@{=>}^<<{\scriptstyle \nu'} (-21,-8); (-19,-12)};
        {\ar@{=>}_<<{\scriptstyle \mu'} (-9,-8); (-11,-12)};
        {\ar@{=>}^<<{} (26,-9); (28,-12)};
        {\ar@{=>}_<<{} (34,-9); (32,-12)};
\endxy
\]
\noindent In fact, this tensor product will be not just a functor, but a homomorphism of bicategories.  This means that it carries some more structure and satisfies some extra axioms, but we will not give these details here.  See the manuscript of Gordon, Power, and Street\;\cite{GPS} for the definition of monoidal bicategory and homomorphism of monoidal bicategory.

\subsection{Enriched Bicategories}\label{enriched}

A monoidal structure, such as the tensor product on $\Span(\Grpd)$ discussed in the previous section, is the crucial ingredient for defining {\em enriched bicategories}.  In particular, given a monoidal bicategory $\mathcal{V}$ with the tensor product $\otimes$, a $\mathcal{V}$-enriched bicategory has for each pair of objects $x,y$, an object $\hom(x,y)$ of $\mathcal{V}$.  Composition involves the tensor product in $\mathcal{{V}}$
\[ \circ\maps\hom(x,y)\otimes \hom(y,z) \to \hom(x,z).\]
\noindent While writing the final draft of this article, we realized that our definition of enriched bicategory is almost identical to one previously given by Carmody\;\cite{Car} and recalled in part by Forcey\;\cite{For}.

After giving the basic structure of enriched bicategories, we state a {\em change of base} theorem, which says which sort of map $f \maps \mathcal{V} \to \mathcal{V}'$ lets us turn a $\mathcal{V}$-enriched bicategory into a $\mathcal{V}'$-enriched bicategory.

Recall that, for each finite group $G$, there is a category of permutation representations $\PermRep(G)$.  Enriched bicategories allow us to define the $\Span(\Grpd)$-enriched bicategory $\Hecke(G)$, which we call the Hecke bicategory.  In particular, the composition by weak pullback of spans fails to satisfy the appropriate coherence conditions of an ordinary bicategory.  The change of base theorem is the main tool employed in proving the our categorification theorem via degroupoidification.

Before giving the definition of an enriched bicategory, we recall the definition of an enriched category --- that is, a category enriched over a monoidal category $\mathcal{V}$\;\cite{Kelly}.  An {\bf enriched category} consists of:
\begin{itemize}
\item a set of objects $x,y,z\ldots$,
\item for each pair of objects $x,y$, an object $\hom(x,y)\in \mathcal{V}$,
\item composition and identity-assigning maps that are morphisms in $\mathcal{V}$.
\end{itemize}
\noindent For example, $\PermRep(G)$ is a category enriched over the monoidal category of vector spaces.

We now define enriched bicategories:

\begin{definition}
Let $\mathcal{V}$ be a monoidal bicategory.  A {\bf $\mathcal{V}$-bicategory} $\mathcal{B}$ consists of the
following data subject to the following axioms:\\
Data:
\begin{itemize}
\item a collection of objects $x,y,z,\ldots$,
\item for every pair of objects $x,y$, a {\bf hom-object} $\hom(x,y)\in\mathcal{V}$,
which we will often denote $(x,y)$,
\item a morphism called {\bf composition}
\[\circ\maps \hom(x,y)\otimes\hom(y,z)\to\hom(x,z)\]
for each triple of objects $x,y,z\in \mathcal{B}$,
\item an {\bf identity-assigning} morphism
\[i_x\maps I\to \hom(x,x)\]
for each object $x\in\mathcal{B}$,
\item an invertible $2$-morphism called the {\bf associator}
\[ \def\objectstyle{\scriptstyle}
  \def\labelstyle{\scriptstyle}
   \xy
   (20,22)*+{(w,x)\otimes ((x,y)\otimes(y,z))}="1";
   (-20,22)*+{((w,x)\otimes (x,y))\otimes(y,z)}="2";
   (-25,0)*+{(w,y)\otimes(y,z)}="4";
   (25,0)*+{(w,x)\otimes(x,z)}="3";
   (0,-20)*+{(w,z)}="5";
        {\ar_{a} "2";"1"};
        {\ar^{1\otimes c} "1";"3"};
        {\ar_{c\otimes 1} "2";"4"};
        {\ar^{c} "3";"5"};
        {\ar_{c} "4";"5"};
        {\ar@2{->}_<<{\scriptstyle \alpha_{wxyz}} (3,6); (-3,3)};
\endxy
\]
for each quadruple of objects $w,x,y,z\in\mathcal{B}$;
\item and invertible $2$-morphisms called the {\bf right unitor} and {\bf left unitor}
\[ \def\objectstyle{\scriptstyle}
  \def\labelstyle{\scriptstyle}
   \xy
		(-10,15)*+{(x,x)\otimes (x,y)}="1";
		(-40,-15)*+{(x,y)}="2";	
		(-10,-15)*+{I\otimes (x,y)}="3";				
		(10,15)*+{(x,y)\otimes (y,y)}="4";
		(40,-15)*+{(x,y)}="5";
		(10,-15)*+{(x,y)\otimes I}="6";
        {\ar@{->}_{c_{xxy}} "1";"2"};
        {\ar@{->}^{r_{xy}} "3";"2"};
        {\ar@{->}_{i_x\otimes 1} "3";"1"};
        {\ar@{->}^{c_{xyy}} "4";"5"};
        {\ar@{->}_{l_{xy}} "6";"5"};
        {\ar@{->}^{1\otimes i_y} "6";"4"};
        {\ar@2{->}_<<{\scriptstyle \rho_{xy}} (-16,-6); (-22,0)};
        {\ar@2{->}_<<{\scriptstyle \lambda_{xy}} (16,-6); (22,0)};
\endxy
\]
for every pair of objects $x,y\in\mathcal{B}$;
\item and axioms given by closed surface diagrams, the more interesting of the two being the {\em permutahedron}\;\;\cite{Car,Hoffnung2}.
\end{itemize}
\end{definition}

Given a monoidal bicategory $\mathcal{V}$, which has only identity $2$-morphisms, then every $\mathcal{V}$-bicategory
is a $\mathcal{V}$-category in the obvious way, and every $\mathcal{V}$-enriched category can be trivially extended to a $\mathcal{V}$-bicategory.  This flexibility will allow us to think of $\PermRep(G)$ as either a $\Vect$-enriched category or as a $\Vect$-enriched bicategory.

Now we state a {\em change of base} construction which allows us to change a $\mathcal{V}$-enriched bicategory to a $\mathcal{V}'$-enriched bicategory.

\begin{claim}\label{changeofbase}
Given a lax-monoidal homomorphism of monoidal bicategories $f\maps \mathcal{V}\to \mathcal{V}'$ and a $\mathcal{V}$-bicategory $\mathcal{B}_\mathcal{V}$, then there is a $\mathcal{V}'$-bicategory
\[ \bar{f}(\mathcal{B}_\mathcal{V}).\]
\end{claim}

A monoidal homomorphism is a map between bicategories preserving the monoidal structure up to isomorphism.\;\cite{GPS,Gurski}.  A {\em lax}-monoidal homomorphism $f$ is a bit more general: it need not preserve the tensor product up to isomorphism.  Instead, it preserves the tensor product only {\em up to a morphism}:
\[ f(x) \otimes' f(y) \to f(x \otimes y),\]
\noindent where the symbol $\otimes'$ is the monoidal product in $\mathcal{V}'$.

The data of the enriched bicategory $\bar{f}(\mathcal{B}_\mathcal{V})$ is straightforward to write down and the proof of the claim is an equally straightforward, yet tedious surface diagram chase.  Here we just point out the most important idea.  The new enriched bicategory $\bar{f}(\mathcal{B}_\mathcal{V})$ has the same objects as $\mathcal{B}_\mathcal{V}$, and for each pair of objects $x,y$, the $\hom$-category of $\bar{f}(\mathcal{B}_\mathcal{V})$ is
\[ \hom_{\bar{f}(\mathcal{B}_\mathcal{V})}(x,y) := f(\hom_{\mathcal{B}_\mathcal{V}}(x,y)).\]

This theorem will allow us to compare the Hecke bicategory, which we define in the next section, to bicategory of spans of finite $G$-sets $\Span(\GSet)$.

\section{A Categorification Theorem}\label{fundamentalsection}

The following sections are devoted to categorifying the Hecke algebroid.  We will show how to obtain the category of permutation representations using the process of {\em degroupoidification}.

\subsection{The Hecke Bicategory}\label{two}

We are now in a position to present the spans of groupoids enriched bicategory $\Hecke(G)$ --- the {\em Hecke bicategory}.

\begin{claim}\label{heckebicategory}
 Given a finite group $G$, there is a $\Span(\Grpd)$-enriched bicategory $\Hecke(G)$ which has:
\begin{itemize}
\item finite $G$-sets $X,Y,Z\ldots$ as objects,
\item for each pair of finite $G$-sets $X, Y$, an object of $\Span$, the action groupoid:
    \[\hom(X,Y) = (X\times Y)\Over G,\]
\item composition
\[ \circ \maps (X\times Y)\Over G \times (Y\times Z)\Over G \to (X\times Z)\Over G\]
is the span of groupoids,
\[
\xymatrix{
& (X\times Y\times Z)\Over G\ar[dr]^{\qquad(p_X\times p_Y)\times (p_Y\times p_Z)}\ar[dl]_{p_X\times p_Z}& \\
(X\times Z)\Over G && (X\times Y)\Over G\times (Y\times Z)\Over G
}
\]
\item for each finite $G$-set $X$, an identity assigning span from the terminal groupoid $1$ to $(X\times X)\Over G$,

\item invertible $2$-morphisms in $\Span(\Grpd)$ assuming the role of the associator and left and right unitors.
\end{itemize}
\end{claim}

Given this structure one needs to check that the axioms of an enriched bicategory are satisfied; however, we will not prove this here.  Combining the degroupoidification functor of Section\;\ref{degroupoidification}, the change of base theorem of Section\;\ref{enriched}, and the enriched bicategory $\Hecke(G)$ described above, we can now state the main theorem.  This is the content of the next section.

\subsection{A Categorification of the Hecke Algebroid}\label{fundamental}

In this section, we describe the relationship between the Hecke algebroid $\PermRep(G)$ of permutation representations of a finite group $G$ and the Hecke bicategory $\Hecke(G)$.  The idea is that for each finite group $G$, the Hecke bicategory $\Hecke(G)$ categorifies $\PermRep(G)$.

We recall the functor {\em degroupoidification}:
\[ \mathcal{D}\maps \Span(\Grpd) \to \Vect\]
\noindent which replaces groupoids with vector spaces and spans of groupoids with linear operators.  With this functor in hand, we can apply the change of base theorem to the $\Span(\Grpd)$-enriched bicategory $\Hecke(G)$.  In other words, for each finite group $G$ there is a $\Vect$-enriched bicategory:
\[\bar{\mathcal{D}}\left(\Hecke(G)\right),\]
\noindent which has
\begin{itemize}
\item permutation representations $X,Y,Z,\ldots$ of $G$ as objects,
\item for each pair of permutation representations $X,Y$, the vector space
\[ \hom(X,Y) = \mathcal{D}\left((X\times Y)\Over G\right)\]
\end{itemize}
\noindent with $G$-orbits of $X\times Y$ as basis.  Of course, a $\Vect$-enriched bicategory is also a $\Vect$-enriched category.  The following is the statement of the main theorem, an equivalence of $\Vect$-enriched categories.

\begin{claim}\label{main}
Given a finite group $G$,
\[\bar{\mathcal{D}}\left(\Hecke(G)\right) \simeq \PermRep(G)\]
as $\Vect$-enriched categories.
\end{claim}

More explicitly, this says that given two permutation representations $\widetilde{X}$ and $\widetilde{Y}$, the vector space of intertwining operators between them can be constructed as the degroupoidification of the groupoid $(X\times Y)\Over G$.

An important corollary of Claim\;\ref{main} is that for certain $G$-sets, which are the flag varieties $X$ associated to Dynkin diagrams, the $\hom$-groupoid $\Hecke(X,X)$ categorifies the associated Hecke algebra.  We will describe these Hecke algebras in Section\;\ref{hecke} and make the relationship to the Hecke bicategory and some of its applications explicit in Section\;\ref{A2}.

\section{Spans of Groupoids and Cocontinuous Functors}\label{alternate}

In the following sections, we sketch the beginning of the project of understanding the relationship between the Hecke bicategory and $\Span(\GSet)$.  For this we will need to introduce the monoidal $2$-category of $\Cocont$ of presheaf categories on groupoids and cocontinuous functors.

\subsection{The Monoidal $2$-Category $\Cocont$}\label{nice}

In Section\;\ref{degroupoidification}, we considered groupoids as categorified vector spaces.  In particular, the isomorphism classes of objects assumed the role of a basis of the corresponding free vector space.  A slightly different point of view, which was discussed at length in\;\cite{HDA7}, assigns to a groupoid the vector space of functions on the set of isomorphism classes of that groupoid.  Thus, promoting functions to functors, we can think of a categorified vector space as the presheaf category on a groupoid.  The objects of the monoidal bicategory described in this section will be defined to be equivalent to such presheaf categories.  We will call such categories {\em nice toposes}.

\begin{definition}
A {\bf presheaf} on a groupoid $\mathcal{G}$ is a contravariant functor from $\mathcal{G}$ to $\Set$.
\end{definition}
\noindent Given a groupoid $\mathcal{G}$, its category of presheaves $\widehat{\mathcal{G}}$ has:
\begin{itemize}
\item presheaves on $\mathcal{G}$ as objects,
\item natural transformations as morphisms.
\end{itemize}

\begin{definition}
A {\bf nice topos} is a category equivalent to the category of presheaves on a (small) groupoid $\mathcal{G}$.
\end{definition}

By the above definition, there is a nice topos $\widehat{\mathcal{G}}$ of presheaves corresponding to any groupoid $\mathcal{G}$.  However, mapping groupoids to these special presheaf categories suggests that nice topoi should have an intrinsic characterization.  To give such a characterization of these topoi we should look to the generalization to topoi of Grothendieck's Galois theory of schemes.  The interested reader is pointed to the survey article\;\cite{Moerdijk} and references therein, although the present paper may be read independently of this survey.  Giving this intrinsic characterization liberates the nice topos from its dependence on a particular groupoid.  In particular, this supports the point of view that nice topoi are the objects of a {\em basis independent} theory of categorified vector spaces.

Following this line of reasoning, the maps between nice topoi are thought of as categorified linear operators.  Thus, they should preserve sums, or more accurately, they should preserve a categorified and generalized notion of `sums' --- colimits.

\begin{definition}
A functor is said to be {\bf cocontinuous} if it preserves all (small) colimits.
\end{definition}

\noindent This suggests that cocontinuous functors might play the role of categorified linear operators.  Indeed, we take such an approach.

In the next section, we will see further support for the analogy: {\em nice topos is to groupoid} as {\em abstract vector space is to vector space with chosen basis} and {\em cocontinuous functor is to span of groupoids} as {\em linear operator is to matrix}.

The monoidal bicategory $\Cocont$ consists of:
\begin{itemize}
\item  nice topoi $\mathcal{D}, \mathcal{E}, \mathcal{F}$ as objects,
\item cocontinuous functors as $1$-morphisms,
\item natural transformations as $2$-morphisms.
\end{itemize}

Objects of $\Cocont$ are categories and the morphisms between them are functors.  Thus, there is a faithful functor from $\Cocont$ to $\Cat$.  It follows that there is a monoidal structure on $\Cocont$ obtained as a restriction of the monoidal structure on $\Cat$.  In particular, the tensor product of nice topoi $\mathcal{E}$ and $\mathcal{F}$ is the cartesian product $\mathcal{E}\times\mathcal{F}$.

In the next section, we will describe the relationship between spans of groupoids and cocontinuous functors between nice topoi.  We use some basic notions of topos theory.

\subsection{From Spans of Groupoids to Cocontinuous Functors}\label{monoidal}

The change of base construction for enriched bicategories offers a new interpretation of the Hecke bicategory.  We have described two closely related monoidal bicategories.  The relationship between groupoids and nice topoi is made manifest as a functor between monoidal bicategories.  Understanding this functor will be the focus of this section.

It is clear from the definition that we can obtain a nice topos by assigning a groupoid $\mathcal{G}$ to its corresponding presheaf category $\widehat{\mathcal{G}}$.  Continuing our analogy with abstract vector spaces and vector spaces with a chosen basis, we will explain how a span of groupoids gives a cocontinuous functor between the corresponding presheaf categories.  In fact, a groupoid can be recovered up to equivalence from its presheaf category, just as a basis can be recovered up to isomorphism from its vector space, but in each case the equivalence or isomorphism is non-canonical.

First, we review some basic ideas from topos theory.  A topos is a category which resembles the category of sets.  Categories of presheaves are examples of {\em Grothendieck topoi}.  In general, a Grothendieck topos is a category of sheaves on a site.  A site is just a category with a notion of a {\em covering} of objects called a {\em Grothendieck topology}\;\cite{Johnstone,MM}.  A familiar example with a particularly simple Grothendieck topology is the category of presheaves on a topological space.

So if a topos is just a special type of category, then how does topos theory differ from category theory?  One answer is that while the morphisms between categories are functors, the morphisms between topoi must satisfy extra properties.  Such a morphism is called a {\em geometric morphism}\;\cite{MM}.  We define the morphisms between nice topoi, although the definition is exactly the same in the more general setting of Grothendieck topoi.

\begin{definition}
A {\bf geometric morphism} $e\maps\mathcal{E}\to\mathcal{F}$ between nice topoi is a pair of functors $e^*\maps \mathcal{F}\to \mathcal{E}$ and $e_*\maps \mathcal{E}\to \mathcal{F}$ such that $e^*$ is left adjoint to $e_*$ and $e^*$ is left exact, i.e., preserves finite limits.  A geometric morphism $e\maps\mathcal{E}\to\mathcal{F}$ is said to be {\bf essential} if there exists a functor $e_!\maps \mathcal{E} \to \mathcal{F}$ which is left adjoint to $e^*$.
\end{definition}

We note a relationship to functors between finite groupoids, which allows us to define cocontinuous functors from spans.  Any functor $f\maps \mathcal{G}\to \mathcal{H}$ defines a geometric morphism between the corresponding presheaf categories:
\[ \widehat{f}\maps \widehat{\mathcal{G}} \to \widehat{\mathcal{H}},\]
\noindent which consists of the functor:
\[ f^*\maps \widehat{\mathcal{H}} \to \widehat{\mathcal{G}},\]
\noindent which pulls presheaves back from $\mathcal{H}$ to $\mathcal{G}$, together with the right adjoint of $f^*$
\[ f_*\maps \widehat{\mathcal{G}} \to \widehat{\mathcal{H}},\]
\noindent which pushes presheaves forward from $\mathcal{G}$ to $\mathcal{H}$.  The particularly important fact is that a geometric morphism induced by a functor between groupoids will always be essential --- that is, there exists a {\em left} adjoint to $f^*$:
\[f_!\maps \widehat{\mathcal{G}} \to \widehat{\mathcal{H}}.\]

\begin{definition}
A {\bf map of geometric morphisms} $\alpha\maps e\Rightarrow f$ is a natural transformation:
\[ \alpha\maps e^*\To f^*.\]
\end{definition}

Using the fact that a functor between groupoids induces an essential geometric morphism, we see that from a span of groupoids:
\[
\xymatrix{
& \mathcal{S}\ar[dl]_q\ar[dr]^p &\\
\mathcal{H} && \mathcal{G}
}
\]
\noindent we can define a functor:
\[q_!p^*\maps \widehat{\mathcal{G}}\to\widehat{\mathcal{H}},\]
\noindent which is a composite of left adjoint functors, and thus, cocontinuous.  Although we do not go into detail at present, using this construction on spans it is not difficult to define a natural transformation between cocontinuous functors from a map of spans.  After checking details, the following claim becomes evident.

\begin{claim}
There is a homomorphism of monoidal bicategories
\[\mathcal{L}\maps \Span \to \Cocont,\]
\noindent which assigns to each groupoid $\mathcal{G}$ its category of sheaves $\widehat{\mathcal{G}}$ and to each span
\[
\xymatrix{
&\mathcal{S}\ar[dl]_{q}\ar[dr]^{p}&\\
\mathcal{H} && \mathcal{G}
}
\]
\noindent the cocontinuous functor $q_!p^*\maps \widehat{G} \to \widehat{H}$.
\end{claim}

In the setting of spans of sets, the homomorphism $\mathcal{L}$ would be analogous to the functor taking a set $X$ to the free vector space with basis $X$, and a span
\[
\xymatrix{
& S\ar[dl]_q\ar[dr]^p &\\
Y && X
}
\]
\noindent to a linear operator between the corresponding vector spaces.

Using the map $\mathcal{L}$, we can apply change of base to the Hecke bicategory to obtain a $\Cocont$-enriched bicategory.  We will discuss the details and benefits of this new structure in the next section.

\subsection{Spans of $G$-Sets as Nice Topoi}\label{presheaf}

In this section we take a closer look at the structure of the $\Cocont$-enriched bicategory $\bar{\mathcal{L}}(\Hecke(G))$.  Our goal is to rephrase the construction as a cocompletion of the bicategory $\Span(\GSet)$ of spans of finite $G$-sets.

\begin{claim}
 Given a finite group $G$, there is a $\Cocont$-enriched bicategory $\bar{\mathcal{L}}(\Hecke(G))$ which has:
\begin{itemize}
\item finite $G$-sets $X,Y,Z\ldots$ as objects,
\item for each pair of finite $G$-sets $X, Y$, an object of $\Nice$,:
    \[\hom(X,Y) = \widehat{(X\times Y)\Over G},\]
\item composition
\[ \circ \maps \widehat{(X\times Y)\Over G} \times \widehat{(Y\times Z)\Over G} \to \widehat{(X\times Z)\Over G}\]
\noindent is the cocontinuous functor $(\pi_{13})_!(\pi_{12}\times \pi_{23})^*$;

\item for each finite $G$-set $X$, an identity assigning cocontinuous functor from the topos $\widehat{1}\simeq \FinSet$ to $\widehat{(X\times X)\Over G}$,

\item invertible $2$-morphisms in $\Cocont$ assuming the role of the associator and left and right unitors.
\end{itemize}
\end{claim}

\noindent The proof of this claim is immediate from the proofs of Claim\;\ref{changeofbase} and Claim\;\ref{heckebicategory}.

It turns out that the $\hom$-categories of the $\Cocont$-enriched bicategory $\bar{\mathcal{L}}(\Hecke(G))$ have a very simple description as spans of $G$-sets and maps of spans.  Given $G$-sets $X$ and $Y$, the product $X\times Y$ is a finite $G$-set, so we can construct the action groupoid $(X\times Y)\Over G$.  The category of presheaves on this groupoid will be a nice topos.

\begin{lemma}\label{grothendieck}
Given a pair of finite $G$-sets $X$, $Y$, the category of spans of $G$-sets and maps of spans of $G$-sets is equivalent to the nice topos $\widehat{(X\times Y)\Over G}$.
\end{lemma}

The construction which proves this lemma is sometimes called the {\em Grothendieck construction}.  Notice that since we consider presheaves valued in $\Set$, rather than only finite sets, we need to allow the apexes of the spans to involve not-necessarily finite sets as well.  The construction says that given a pair of finite $G$-sets $X$ and $Y$, presheaves on $(X\times Y)\Over G$ are spans from $X$ to $Y$ and natural transformations are maps of spans.  We sketch the proof of this lemma now.

\begin{proof}(Sketch)
Given a span of finite $G$-sets:
\[
\xymatrix{
& S\ar[dl]_q\ar[dr]^{p} &\\
Y && X
}
\]
\noindent there is a presheaf on $(X\times Y)\Over G$, which we can think of approximately as a categorified matrix of natural numbers, i.e., a matrix of sets.  Each object $(x,y)$ determines an entry in the matrix, and the entries are the sets $S_{x,y} = p^{-1}(x)\cap q^{-1}(y)$ defined in Section\;\ref{matrix}.  For each morphism $(g,(x,y))\maps (x,y)\to (x',y')$, we define a function from $S_{x',y'}$ to $S_{x,y}$ by the action of $g^{-1}$ on the $G$-set $S_{x',y'}$.  Thus, we get a presheaf from the span $S$.

Now from a map of spans of finite $G$-sets:
\[
\xymatrix{
& S\ar[dr]^p\ar[dl]_q\ar[dd]^f&\\
Y&&X\\
& S'\ar[ul]^{q'}\ar[ur]_{p'} &
}
\]
\noindent we construct a natural transformation between the presheaves corresponding to $S$ and $S'$.  For each object $(x,y)$ of $(X\times Y)\Over G$, the component of the natural transformation takes an element $s \in S_{x,y}$ to $f(s) \in S'_{x,y}$.  Since $f$ is $G$-equivariant, the naturality squares commute.

It is not difficult to check that this process defines an equivalence of categories.  We only need to build a functor in the other direction such that the respective composites are naturally isomorphic to the identity.  We construct such a functor here.

Consider a presheaf on the groupoid $(X\times Y)\Over G$, i.e., an assignment of a set to each $G$-orbit in $X\times Y$.  Projecting out of the product, we obtain a span of $G$-sets by assigning the set over a pair $(x,y)$ the value of the presheaf on the $G$-orbit of this pair.  Next, consider a natural transformation from a presheaf $P$ to a presheaf $Q$.  We define a map between the corresponding spans denoted by their apex sets $S_P$ and $S_Q$.   Each element $p\in S_P$ is assigned by the span a pair $(x,y)$, and to this pair, the natural transformation assigns a function $f \maps P(x,y) \to Q(x,y)$.  Since $p\in P(x,y)\subseteq S_P$ and $f(p)\in Q(x,y)\subseteq S_Q$, we have defined a map of spans.
\end{proof}

We have seen that the $\hom$-categories of the $\Cocont$-enriched bicategory $\bar{\mathcal{L}}(\Hecke(G))$ are actually categories whose objects are spans of $G$-sets and whose morphisms are maps between such spans.  In particular, this bicategory is the cocompletion of $\Span(\GSet)$.  In describing categorified Hecke algebras we will present a groupoid and a span of groupoids.  It will be useful when we consider a specific example to consider the presheaf category in place of the groupoid and to think of each presheaf as a span of $G$-sets.

\section{The Categorified Hecke Algebra and Zamolodchikov Equation}\label{applications}
The main claim of this paper is the existence of a categorification of the Hecke algebroid by the Hecke bicategory $\Hecke(G)$.  This is, in fact, a statement about categorified Hecke algebras.  In attempting to make connections between the categorified Hecke algebra and knot theory, the bicategory $\Span(\GSet)$ allows a hands-on approach to the categorified generators and the isomorphisms of spans arising from the defining equations of the Hecke algebra.  We show that in certain cases these are Yang-Baxter operators that satisfy the Zamolodchikov tetrahedron equation. The hope is that these Yang-Baxter equations will lead to interesting braided monoidal $2$-categories\;\cite{KV}.

\subsection{The Hecke Algebra}\label{hecke}

There are several well-known equivalent descriptions of the Hecke algebra $\mathcal{H}(\Gamma,q)$ obtained from a Dynkin diagram $\Gamma$ and a prime power $q$.  One description of the Hecke algebra is a $q$-deformation of the group algebra of the Coxeter group of $\Gamma$.  A standard example of a Coxeter group associated to a Dynkin diagram is the symmetric group on $n$ letters $S_n$, which is the Coxeter group of the $A_{n-1}$ Dynkin diagram.  We will return to this definition in Section\;\ref{A2} and see that it lends itself to combinatorial applications of the Hecke algebra.  This combinatorial aspect comes from the close link between the Coxeter group and its associated Coxeter complex, a finite simplicial complex that plays an essential role in the theory of buildings\;\cite{Brown}.

Hecke algebras have an alternative definition as algebras of intertwining operators between certain coinduced representations\;\cite{Bump}.  Given a Dynkin diagram $\Gamma$ and prime power $q$, there is an associated simple algebraic group $G = G(\Gamma,q)$.  Choosing a Borel subgroup $B \subset G$, i.e., a maximal solvable subgroup, we can construct the corresponding flag complex $X = G/B$, a transitive $G$-set.

Now, for a finite group $G$ and a representation $V$ of a subgroup $H\subset G$, the {\em coinduced representation} of $G$ from $H$ is defined as the $V$-valued functions on $G$, which commute with the action of $H$:

\[{\rm CoInd}^G_H = \{f\maps G\to V\mid h\cdot f(x) = f(hx)\}\]

\noindent The action of $h\in G$ is defined on a function $f\maps G\to V$ as $h\cdot f(x) = f(xh^{-1})$.  A standard fact about finite groups says that the representation coinduced from the trivial representation of any subgroup is isomorphic to the permutation representation on the cosets of that subgroup.

Thus, from the trivial representation of a Borel subgroup $B$, we obtain the permutation representation $\widetilde{X}$ on the cosets of $B$, i.e., the flag complex $X$.  Then the Hecke algebra is defined as the algebra of intertwining operators from $\widetilde{X}$ to itself:
\[ \mathcal{H}(\Gamma,q) := \PermRep(G)(\widetilde{X},\widetilde{X}),\]
\noindent where $G = G(\Gamma,q)$ and we use the notation $\mathcal{C}(A,B)$ to denote $\hom(A,B)$ in the category $\mathcal{C}$.

Given this definition of the Hecke algebra, we have an immediate corollary to Claim\;\ref{main}:

\begin{claim}\label{maincorollary}
Given a Dynkin diagram $\Gamma$ and prime power $q$, denote $G = G(\Gamma,q)$.  Then the  $\hom$-category $\Hecke(G)(X,X)$ of the Hecke bicategory categorifies the Hecke algebra $\mathcal{H}(\Gamma,q)$.
\end{claim}

\subsection{Categorified Generators and the Zamolodchikov Equation}\label{A2}

Now that we have seen a categorification of Hecke algebras abstractly as a corollary, we look at a concrete example.  The categorified Hecke algebra is particularly easy to understand as living inside the bicategory $\Span(\GSet)$.

While we found it useful earlier to view Hecke algebras as algebras of intertwining operators, viewing the Hecke algebra as a $q$-deformation of a Coxeter group\;\cite{Humphreys} is helpful in examples.

Any Dynkin diagram gives rise to a simple Lie group, and the Weyl group of this simple Lie group is a Coxeter group.  Let $\Gamma$ be a Dynkin diagram.  We write $d \in \Gamma$ to mean that $d$ is a
dot in this diagram.  Associated to each unordered pair of dots $d,d'
\in \Gamma$ is a number $m_{dd'} \in \{2,3,4,6\}$.  In the usual Dynkin
diagram conventions:
\begin{itemize}
\item
$m_{dd'} = 2$ is drawn as no edge at all,
\item
$m_{dd'} = 3$ is drawn as a single edge,
\item
$m_{dd'} = 4$ is drawn as a double edge,
\item
$m_{dd'} = 6$ is drawn as a triple edge.
\end{itemize}

For any prime power $q$, our Dynkin diagram $\Gamma$ defines a Hecke algebra.  The {\em Hecke algebra} $\mathcal{H}(\Gamma,q)$ corresponding to this data is the associative algebra with one generator $\sigma_d$ for each $d \in \Gamma$,
and relations:
\[\sigma_d^2 = (q-1)\sigma_d + q \]
for all $d \in \Gamma$, and
\[\sigma_d\sigma_{d'}\sigma_d\cdots
= \sigma_{d'} \sigma_d \sigma_{d'}\cdots \]
for all $d,d'\in \Gamma$, where each side has $m_{dd'}$ factors.

When $q = 1$, this Hecke algebra is simply the group algebra of
the Coxeter group associated to $\Gamma$ --- that is, the group
with one generator $s_d$ for each dot $d \in \Gamma$, and relations
\[ s_d^2 = 1,  \qquad (s_d s_{d'})^{m_{dd'}} = 1 .\]
So, the Hecke algebra can be thought of as a $q$-deformation of this
Coxeter group.

We recall the flag complex $X = G/B$ from Section\;\ref{hecke} is a finite set equipped with a transitive action of the finite group $G$.  Starting from just this $G$-set $X$, we can see an explicit picture of the categorified Hecke algebra of spans of $G$-sets from $X$ to $X$.

The key is that for each dot $d \in \Gamma$ there is a special
span of $G$-sets that corresponds to the generator $\sigma_d \in \mathcal{H}(\Gamma,q)$.  To illustrate these ideas, let us consider the simplest nontrivial example, the Dynkin diagram $A_2$:
\[\xymatrix{ \bullet\ar@{-}[r] & \bullet}\]
The Hecke algebra associated to $A_2$ has two generators, which we
call $P$ and $L$, for reasons soon to be revealed:
\[        P = \sigma_1 , \qquad L = \sigma_2 . \]
\noindent To make the connection to the description of the Hecke algebra as an algebra of intertwining operators explicit, we can choose a basis of $\widetilde{X}$ corresponding to the elements of the Coxeter group $S_3$ and write these generators as the following matrices, or intertwining operators.
\[ P = \left( \begin{array}{cccccc}
0 & q & 0 & 0 & 0 & 0 \\
1 & q-1 & 0 & 0 & 0 & 0 \\
0 & 0 & 0 & q & 0 & 0 \\
0 & 0 & 1 & q-1 & 0 & 0\\
0 & 0 & 0 & 0 & q-1 & 1\\
0 & 0 & 0 & 0 & q & 0
\end{array} \right)\] 
\[ L = \left( \begin{array}{cccccc}
0 & 0 & 0  & 0 & 0 & q\\
0 & 0 & q  & 0 & 0 & 0\\
0 & 1 & q-1  & 0 & 0 & 0\\
0 & 0 & 0  & q-1 & 1 & 0\\
0 & 0 & 0  & q & 0 & 0\\
1 & 0 & 0  & 0 & 0 & q-1
\end{array} \right)\] 
The relations are
\[     P^2 = (q-1) P + q, \qquad L^2= (q-1)L + q,  \qquad PLP=LPL  .\]
It follows that this Hecke algebra is a quotient of the group algebra
of the 3-strand braid group, which has two generators $P$ and $L$, which we can draw as braids in $3$-dimensional space:
\[
\begin{picture}(60,70)(60,0)
 \includegraphics[scale=0.5]{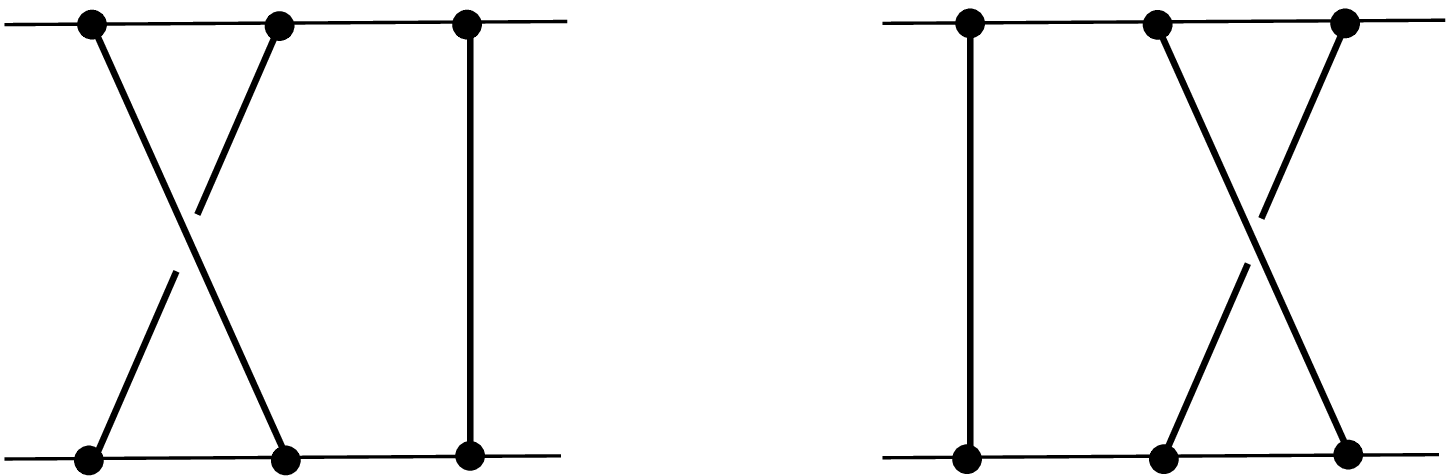}
\put(-235,30){$P = $}
\put(-110,30){$L = $}
\end{picture}
\]
\noindent and one relation $PLP = LPL$:
\[
\begin{picture}(45,90)(45,0)
 \includegraphics[scale=0.5]{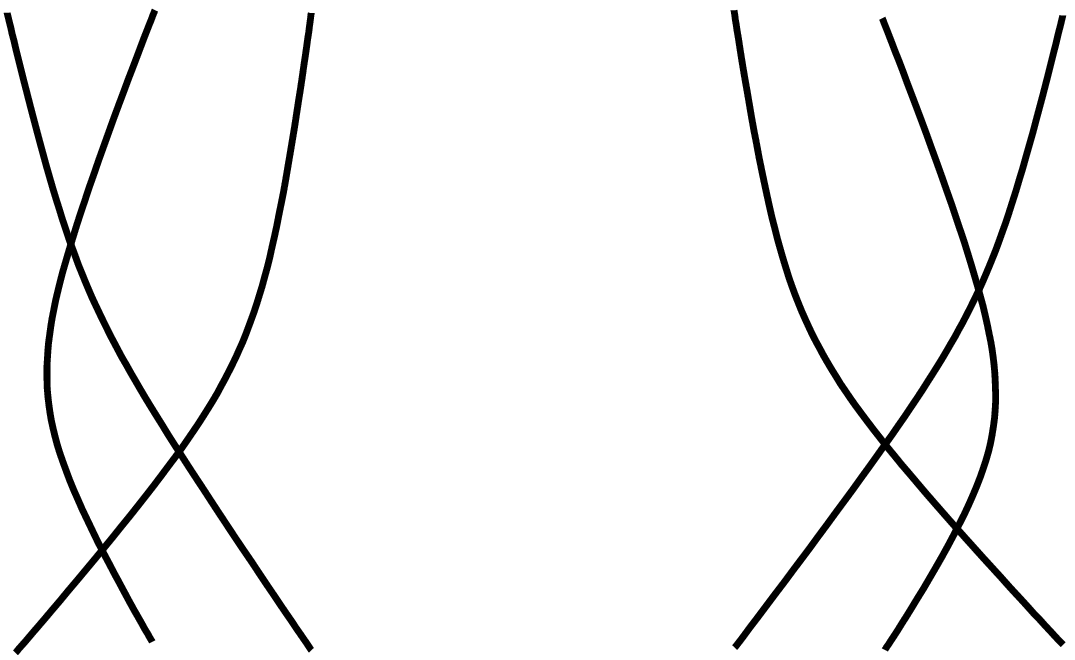}
  \put(-80,45){$=$}
\end{picture}
\]
\noindent called the {\em Yang--Baxter equation} or
{\em third Reidemeister move}.  This is why Jones could use traces on
the $A_n$ Hecke algebras to construct invariants of knots\;\cite{Jones}.  In light of the success of Khovanov homology --- a categorification of the Jones polynomial --- this connection to knot theory makes it especially interesting to categorify Hecke algebras.

So, let us see what the categorified Hecke algebra looks like, and
where the Yang--Baxter equation comes from.  The algebraic group
corresponding to the $A_2$ Dynkin diagram and the prime power $q$ is
$G = \SL(3,\F_q)$, and we can choose the Borel subgroup $B$ to consist
of upper triangular matrices in $\SL(3,\F_q)$.  Recall that a complete flag in the vector space $\F_q^3$ is a pair of subspaces
\[           0 \subset V_1 \subset V_2 \subset \F_q^3  . \]
The subspace $V_1$ must have dimension one, while $V_2$ must have
dimension two.  Since $G$ acts transitively on the set of complete
flags and $B$ is the subgroup stabilizing a chosen flag, the flag
variety $X = G/B$ in this example is just the set of complete flags in
$\F_q^3$---hence its name.

We can think of $V_1 \subset \F_q^3$ as a point in the projective
plane $\F_q{\mathrm P}^2$, and $V_2 \subset \F_q^3$ as a line in this
projective plane.  From this viewpoint, a complete flag is a chosen
point lying on a chosen line in $\F_q {\mathrm P}^2$.  This viewpoint
is natural in the theory of `buildings', where each Dynkin diagram
corresponds to a type of geometry\;\cite{Brown}.  Each dot in
the Dynkin diagram then stands for a `type of geometrical figure',
while each edge stands for an `incidence relation'.  The $A_2$ Dynkin
diagram corresponds to projective plane geometry.  The dots in this
diagram stand for the figures `point' and `line':
\[\xymatrix{ {\rm point} \; \bullet \ar@{-}[r] &
\bullet \; {\rm line} }\]
The edge in this diagram stands for the incidence
relation `the point $p$ lies on the line $\ell$'.

We can think of $P$ and $L$ as special elements of the $A_2$ Hecke
algebra, as already described.  But when we categorify the Hecke
algebra, $P$ and $L$ correspond to irreducible spans of $G$-sets -- that is, spans that are not coproducts of two non-trivial spans of $G$-sets.  Let us describe these spans and explain how the Hecke algebra relations arise in this categorified setting.

The objects $P$ and $L$ can be defined by giving irreducible spans
of $G$-sets:
\[
\xymatrix{
    & P \ar[dl]\ar[dr] &    &&     & L\ar[dl]\ar[dr] & \\
X &                  &X && X &                 &X
}
\]

In general, any span of $G$-sets
\[
\xymatrix{
& S \ar[dl]_{q} \ar[dr]^{p} & \\
X & & X}
\]
such that $q \times p \maps S \to X \times X$ is injective can be
thought of as $G$-invariant binary relation between elements of $X$.
Irreducible $G$-invariant spans are always injective in this sense.
So, such spans can also be thought of as $G$-invariant relations between flags.  In these terms, we define $P$ to be the relation that says two flags have the same line, but different points:
\[ P = \{((p,\ell),(p',\ell)) \in X \times X \mid p\neq p'\}
\]
Similarly, we think of $L$ as a relation saying two flags
have different lines, but the same point:
\[ L =
\{((p,\ell),(p,\ell')) \in X \times X \mid \ell\neq \ell'\}. \]
Given this, we will check:
\[ P^2 \cong (q-1) \times P + q \times 1, \qquad L^2 \cong (q-1) \times L +
q \times 1, \qquad PLP \cong LPL . \]
Here both sides refer to spans of $G$-sets.  Addition of spans is defined using the coproduct, while $1$ denotes the identity span from $X$ to $X$.  We use `$q$' to stand for a fixed $q$-element set, and similarly for `$q-1$'.  We compose spans of $G$-sets using the ordinary pullback.

To check the existence of the first two isomorphisms above, we just
need to count.  In $\mathbb{F}_q\mathrm{P}^2$, the are $q+1$ points on any line.  So, given a flag we can change the point in $q$ different ways.  To change it again, we have a choice: we can either send it back to the original point, or change it to one of the $q-1$ other points.
So, $P^2 \cong (q-1) \times P + q \times 1$.  Since there are also
$q+1$ lines through any point, similar reasoning shows that $L^2 \cong
(q-1) \times L + q \times 1$.

The Yang--Baxter isomorphism
\[          PLP \cong LPL  \]
is more interesting.  We construct it as follows.  First consider the
left-hand side, $PLP$.  So, start with a complete flag called
$(p_1,\ell_1)$:
\hfill \break
\begin{center}
\begin{picture}(50,25)
  \includegraphics[scale=0.35]{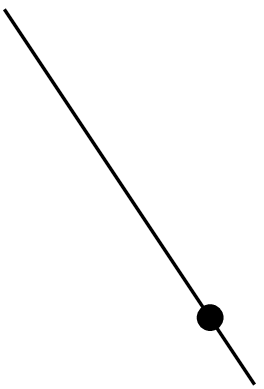}
  \put(-20,2){$p_1$}
  \put(0,-10){$\ell_1$}
\end{picture}
\end{center}
\medskip
Then, change the point to obtain a flag $(p_2,\ell_1)$.
Next, change the line to obtain a flag $(p_2,\ell_2)$.
Finally, change the point once more, which gives us the flag $(p_3,\ell_2)$:
\begin{center}
  \begin{picture}(250,45)
  \includegraphics[scale=0.35]{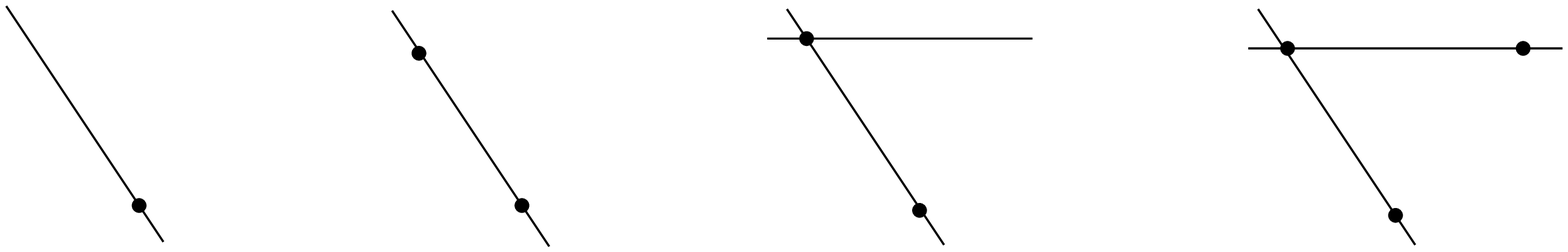}
  \put(-245,4){$p_1$}
  \put(-225,-10){$\ell_1$}
  \put(-185,4){$p_1$}
  \put(-165,-10){$\ell_1$}
  \put(-200,25){$p_2$}
  \put(-120,4){$p_1$}
  \put(-100,-10){$\ell_1$}
  \put(-130,23){$p_2$}
  \put(-80,35){$\ell_2$}
  \put(-40,4){$p_1$}
  \put(-20,-10){$\ell_1$}
  \put(-55,23){$p_2$}
  \put(5,35){$\ell_2$}
  \put(-15,40){$p_3$}
  \end{picture}
\end{center}
\medskip
\noindent
The figure on the far right is a typical element of $PLP$:
\[ ((p_1,\ell_1),(p_2,\ell_1),(p_2,\ell_2),(p_3,\ell_2)) \textrm{ such that } p_1\neq p_2, p_2\neq p_3, \ell_1\neq\ell_2.\]

On the other hand, consider $LPL$.  So, start with the same flag as
before, but now change the line, obtaining $(p_1,\ell'_2)$.  Next
change the point, obtaining the flag $(p'_2,\ell'_2)$.  Finally,
change the line once more, obtaining the flag $(p'_2,\ell'_3)$:
\hfill \break
\medskip
\begin{center}
\begin{picture}(240,40)
  \includegraphics[scale=0.35]{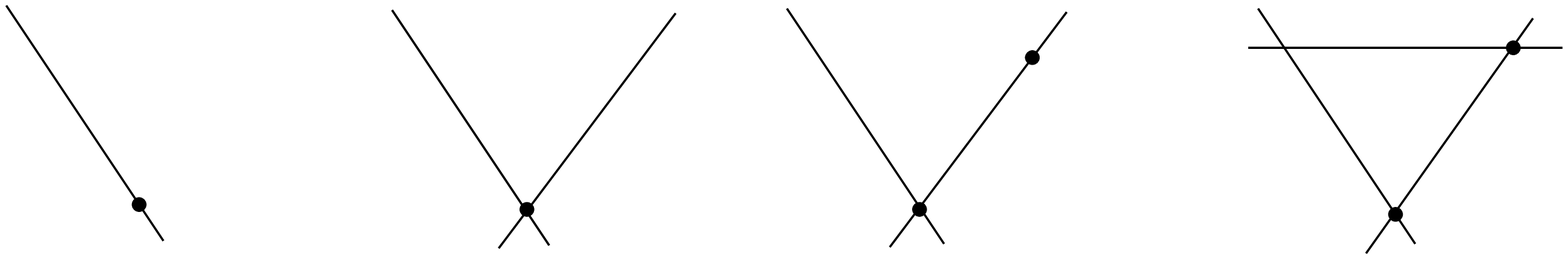}
  \put(-245,3){$p_1$}
  \put(-225,-10){$\ell_1$}
  \put(-185,3){$p_1$}
  \put(-170,-10){$\ell_1$}
  \put(-150,45){$\ell_2'$}
  \put(-120,3){$p_1$}
  \put(-100,-10){$\ell_1$}
  \put(-85,45){$\ell_2'$}
  \put(-85,20){$p_2'$}
  \put(-45,4){$p_1$}
  \put(-25,-10){$\ell_1$}
  \put(-10,43){$\ell_2'$}
  \put(-10,20){$p_2'$}
  \put(-60,25){$\ell_3'$}
\end{picture}
\end{center}
\medskip
\noindent
The figure on the far right is a typical element of $LPL$.

Now, the axioms of projective plane geometry say that any two distinct
points lie on a unique line, and any two distinct lines intersect in a
unique point.  So, any figure of the sort shown on the left below
determines a unique figure of the sort shown on the right, and vice
versa:
\hfill\break
\begin{center}
\begin{picture}(150,50)
\includegraphics[scale=0.50]{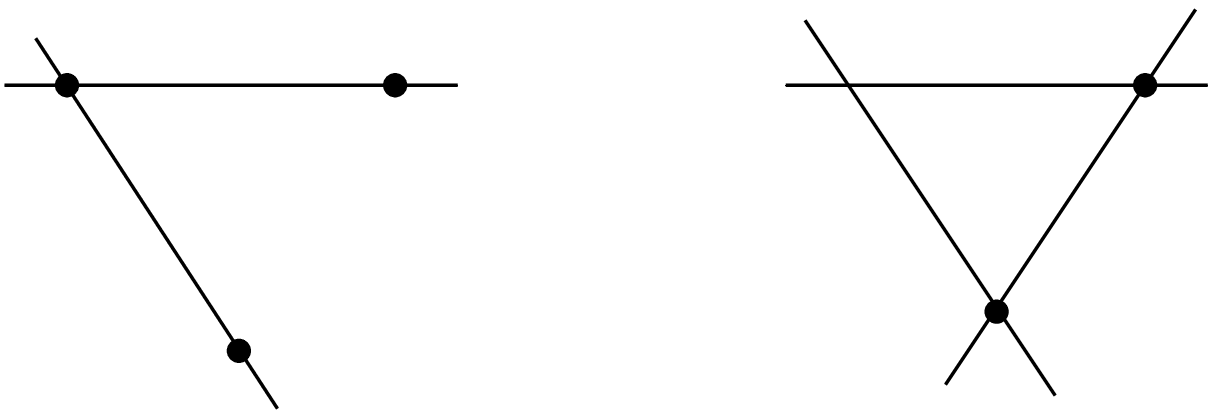}
\end{picture}
\end{center}
\medskip
\noindent
Comparing this with the pictures above, we see this bijection induces
an isomorphism of spans $PLP \cong LPL$.  So, we have derived the
Yang--Baxter isomorphism from the axioms of projective plane geometry!

The above discussion helps illuminate the occurrence of the Yang--Baxter {\em equation} in the generators and relations description of the Hecke algebra.  We have seen that the categorified setting allows us to view these equations as {\em isomorphisms} of spans of $G$-sets.  As such, these {\em Yang-Baxter operators} satisfy an equation of their own -- the {\em Zamolodchikov tetrahedron equation}\;\cite{KV}.  However, this equation only appears in the categorified $A_n$ Hecke algebra, for $n\geq 3$.  We can assign braids on four strands to the generators of the $A_3$ Hecke algebra:
\[
\begin{picture}(110,70)(110,0)
 \includegraphics[scale=0.5]{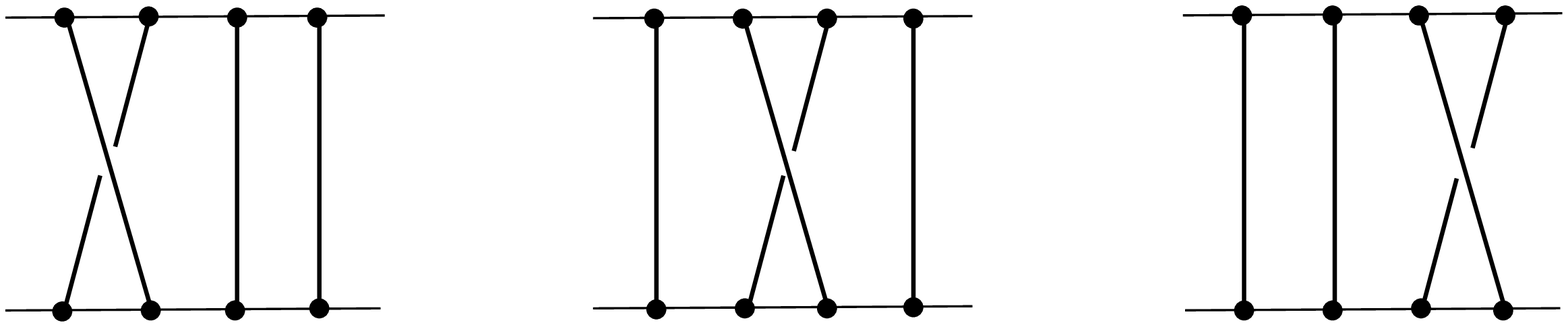}
\put(-340,30){$P = $}
\put(-235,30){$L = $}
\put(-110,30){$S = $}
\end{picture}
\]
\noindent where composition of spans, or multiplication in the Hecke algebra, corresponds to stacking of braid diagrams.  Then we can express the Zamolodchikov equation -- as an equation in the categorified Hecke algebra -- in the form of a commutative diagram of braids\;\cite{BaezCrans,CarterSaito}:

\[
\xy 0;/r.11pc/: 
(0,50)*+{
\xy 
(-15,-20)*{}="T1";
(-5,-20)*{}="T2";
(5,-20)*{}="T3";
(15,-20)*{}="T4";
(-14,20)*{}="B1";
(-5,20)*{}="B2";
(5,20)*{}="B3";
(15,20)*{}="B4";
"T1"; "B4" **\crv{(-15,-7) & (15,-5)}
\POS?(.25)*{\hole}="2x" \POS?(.47)*{\hole}="2y" \POS?(.6)*{
\hole}="2z";
"T2";"2x" **\crv{(-4,-12)};
"T3";"2y" **\crv{(5,-10)};
"T4";"2z" **\crv{(16,-9)};
(-15,-5)*{}="3x";
"2x"; "3x" **\crv{(-18,-10)};
"3x"; "B3" **\crv{(-13,0) & (4,10)}
\POS?(.3)*{\hole}="4x" \POS?(.53)*{\hole}="4y";
"2y"; "4x" **\crv{};
"2z"; "4y" **\crv{};
(-15,10)*{}="5x";
"4x";"5x" **\crv{(-17,6)};
"5x";"B2" **\crv{(-14,12)}
\POS?(.6)*{\hole}="6x";
"6x";"B1" **\crv{(-14,18)};
"4y";"6x" **\crv{(-8,10)};
\endxy
}="T";
(-40,30)*+{
\xy 
(-15,-20)*{}="b1";
(-5,-20)*{}="b2";
(5,-20)*{}="b3";
(14,-20)*{}="b4";
(-14,20)*{}="T1";
(-5,20)*{}="T2";
(5,20)*{}="T3";
(15,20)*{}="T4";
"b1"; "T4" **\crv{(-15,-7) & (15,-5)}
\POS?(.25)*{\hole}="2x" \POS?(.47)*{\hole}="2y" \POS?(.65)*{ \hole}="2z";
"b2";"2x" **\crv{(-5,-15)};
"b3";"2y" **\crv{(5,-10)};
"b4";"2z" **\crv{(14,-9)};
(-15,-5)*{}="3x";
"2x"; "3x" **\crv{(-15,-10)};
"3x"; "T3" **\crv{(-15,15) & (5,10)}
\POS?(.38)*{\hole}="4y" \POS?(.65)*{\hole}="4z";
"T1";"4y" **\crv{(-14,16)};
"T2";"4z" **\crv{(-5,16)};
"2y";"4z" **\crv{(-10,3) & (10,2)} \POS?(.6)*{\hole}="5z";
"4y";"5z" **\crv{(-5,5)};
"5z";"2z" **\crv{(5,4)};
\endxy }="TL";
(-75,0)*+{
\xy 
(-14,20)*{}="T1";
(-4,20)*{}="T2";
(4,20)*{}="T3";
(15,20)*{}="T4";
(-15,-20)*{}="B1";
(-5,-20)*{}="B2";
(5,-20)*{}="B3";
(15,-20)*{}="B4";
"B1";"T4" **\crv{(-15,5) & (15,-5)}
\POS?(.25)*{\hole}="2x" \POS?(.49)*{\hole}="2y" \POS?(.65)*{ \hole}="2z";
"2x";"T3" **\crv{(-20,10) & (5,10) }
\POS?(.45)*{\hole}="3y" \POS?(.7)*{\hole}="3z";
"2x";"B2" **\crv{(-5,-14)};
"T1";"3y" **\crv{(-16,17)};
"T2";"3z" **\crv{(-5,17)};
"3z";"2z" **\crv{};
"3y";"2y" **\crv{};
"B3";"2z" **\crv{ (5,-5) &(20,-10)}
\POS?(.4)*{\hole}="4z";
"2y";"4z" **\crv{(6,-8)};
"4z";"B4" **\crv{(15,-15)};
\endxy }="ML";
(-40,-30)*+{
\xy 
(-14,20)*{}="T1";
(-4,20)*{}="T2";
(4,20)*{}="T3";
(15,20)*{}="T4";
(-15,-20)*{}="B1";
(-5,-20)*{}="B2";
(5,-20)*{}="B3";
(15,-20)*{}="B4";
"B1";"T4" **\crv{(-15,-5) & (15,5)}
\POS?(.38)*{\hole}="2x" \POS?(.53)*{\hole}="2y" \POS?(.7)*{\hole}="2z";
"T1";"2x" **\crv{(-15,5)};
"2y";"B2" **\crv{(10,-10) & (-6,-10)} \POS?(.45)*{\hole}="4x";
"2z";"B3" **\crv{ (15,0)&(15,-10) & (6,-16)} \POS?(.7)*{\hole}="5x";
"T3";"2y" **\crv{(5,10)& (-6,18) }
\POS?(.5)*{\hole}="3x";
"T2";"3x" **\crv{(-5,15)};
"3x";"2z" **\crv{(7,11)};
"2x";"4x" **\crv{(-3,-7)};
"4x";"5x" **\crv{};
"5x";"B4" **\crv{(15,-15)};
\endxy }="BL";
(40,30)*+{
\xy 
(14,-20)*{}="T1";
(4,-20)*{}="T2";
(-4,-20)*{}="T3";
(-15,-20)*{}="T4";
(15,20)*{}="B1";
(5,20)*{}="B2";
(-5,20)*{}="B3";
(-15,20)*{}="B4";
"B1";"T4" **\crv{(15,5) & (-15,-5)}
\POS?(.38)*{\hole}="2x" \POS?(.53)*{\hole}="2y" \POS?(.7)*{\hole}="2z";
"T1";"2x" **\crv{(15,-5)};
"2y";"B2" **\crv{(-10,10) & (6,10)}
\POS?(.45)*{\hole}="4x";
"2z";"B3" **\crv{ (-15,0)&(-15,10) & (-6,16)}
\POS?(.7)*{\hole}="5x";
"T3";"2y" **\crv{(-5,-10)& (6,-18) }
\POS?(.5)*{\hole}="3x";
"T2";"3x" **\crv{(5,-15)};
"3x";"2z" **\crv{(-7,-11)};
"2x";"4x" **\crv{(3,7)};
"4x";"5x" **\crv{};
"5x";"B4" **\crv{(-15,15)};
\endxy }="TR";
(75,0)*+{
\xy 
(14,-20)*{}="T1";
(4,-20)*{}="T2";
(-4,-20)*{}="T3";
(-15,-20)*{}="T4";
(15,20)*{}="B1";
(5,20)*{}="B2";
(-5,20)*{}="B3";
(-15,20)*{}="B4";
"B1";"T4" **\crv{(15,-5) & (-15,5)}
\POS?(.25)*{\hole}="2x" \POS?(.49)*{\hole}="2y" \POS?(.65)*{ \hole}="2z";
"2x";"T3" **\crv{(20,-10) & (-5,-10) }
\POS?(.45)*{\hole}="3y" \POS?(.7)*{\hole}="3z";
"2x";"B2" **\crv{(5,14)};
"T1";"3y" **\crv{(16,-17)};
"T2";"3z" **\crv{(5,-17)};
"3z";"2z" **\crv{};
"3y";"2y" **\crv{};
"B3";"2z" **\crv{ (-5,5) &(-20,10)}
\POS?(.4)*{\hole}="4z";
"2y";"4z" **\crv{(-6,8)};
"4z";"B4" **\crv{(-15,15)};
\endxy }="MR";
(40,-30)*+{
\xy 
(15,20)*{}="b1";
(5,20)*{}="b2";
(-5,20)*{}="b3";
(-14,20)*{}="b4";
(14,-20)*{}="T1";
(5,-20)*{}="T2";
(-5,-20)*{}="T3";
(-15,-20)*{}="T4";
"b1"; "T4" **\crv{(15,7) & (-15,5)}
\POS?(.25)*{\hole}="2x" \POS?(.47)*{\hole}="2y" \POS?(.65)*{ \hole}="2z";
"b2";"2x" **\crv{(5,15)};
"b3";"2y" **\crv{(-5,10)};
"b4";"2z" **\crv{(-14,9)};
(15,5)*{}="3x";
"2x"; "3x" **\crv{(15,10)};
"3x"; "T3" **\crv{(15,-15) & (-5,-10)}
\POS?(.38)*{\hole}="4y" \POS?(.65)*{\hole}="4z";
"T1";"4y" **\crv{(14,-16)};
"T2";"4z" **\crv{(5,-16)};
"2y";"4z" **\crv{(10,-3) & (-10,-2)} \POS?(.6)*{\hole}="5z";
"4y";"5z" **\crv{(5,-5)};
"5z";"2z" **\crv{(-5,-4)};
\endxy }="BR";
(0,-50)*+{
\xy 
(15,20)*{}="T1";
(5,20)*{}="T2";
(-5,20)*{}="T3";
(-15,20)*{}="T4";
(15,-20)*{}="B1";
(5,-20)*{}="B2";
(-5,-20)*{}="B3";
(-15,-20)*{}="B4";
"T1"; "B4" **\crv{(15,7) & (-15,5)}
\POS?(.25)*{\hole}="2x" \POS?(.45)*{\hole}="2y" \POS?(.6)*{\hole}="2z";
"T2";"2x" **\crv{(4,12)};
"T3";"2y" **\crv{(-5,10)};
"T4";"2z" **\crv{(-16,9)};
(15,5)*{}="3x";
"2x"; "3x" **\crv{(18,10)};
"3x"; "B3" **\crv{(13,0) & (-4,-10)}
\POS?(.3)*{\hole}="4x" \POS?(.53)*{\hole}="4y";
"2y"; "4x" **\crv{};
"2z"; "4y" **\crv{};
(15,-10)*{}="5x";
"4x";"5x" **\crv{(17,-6)};
"5x";"B2" **\crv{(14,-12)}
\POS?(.6)*{\hole}="6x";
"6x";"B1" **\crv{};
"4y";"6x" **\crv{};
\endxy }="B";
(-20,65)*{}="X1";
(-35,55)*{}="X2"; 
{\ar@{=>} "X1";"X2"}; 
(20,65)*{}="X1"; 
(35,55)*{}="X2"; 
{\ar@{=>} "X1";"X2"}; 
(60,40)*{}="X1"; 
(75,25)*{}="X2"; 
{\ar@{=>} "X1";"X2"}; 
(-60,40)*{}="X1"; 
(-75,25)*{}="X2"; 
{\ar@{=>} "X1";"X2"};
(-60,-40)*{}="X2";
(-75,-25)*{}="X1";
{\ar@{=>} "X1";"X2"};
(60,-40)*{}="X2";
(75,-25)*{}="X1";
{\ar@{=>} "X1";"X2"};
(-20,-65)*{}="X2";
(-35,-55)*{}="X1";
{\ar@{=>} "X1";"X2"};
(20,-65)*{}="X2";
(35,-55)*{}="X1";
{\ar@{=>} "X1";"X2"};
\endxy \]

This is just the beginning of a wonderful story involving Dynkin diagrams of more general types, incidence geometries, logic, braided monoidal $2$-categories\;\cite{BaezNeuchl,Mc}, knot invariants, topological quantum field theories, geometric representation theory, and more!

\section*{Acknowledgements}

A tremendous debt is owed to John Baez, James Dolan and Todd Trimble for initiating an interesting and creative long-term endeavour, as well as for generously sharing their time and ideas.  Also to the contributors to the $n$-Category Caf\'e, especially John Baez, Bruce Bartlett, Denis-Charles Cisinski, Tom Leinster, Mike Shulman, Todd Trimble, Urs Schreiber and Simon Willerton for taking an interest and sharing and clarifying so many ideas.  To David Ben-Zvi, Anthony Licata, Urs Schreiber and the members of the {\em Journal Club} for introducing the author to the larger world of geometric function theories.  Thanks to Nick Gurski for patiently explaining the basics of tricategories; Mikhail Khovanov and Aaron Lauda for generously hosting me at Columbia University while many of these ideas were developed; Jim Stasheff for encouragement on a very early draft; Julia Bergner and Christopher Walker for patiently listening and contributing to many useful conversations; Rick Blute, Alistair Savage, Pieter Hofstra, Phil Scott and others for patiently listening and commenting week after week, while a number of ideas were worked out in our seminar at the University of Ottawa; and to Dorette Pronk for her invitation to Dalhousie University, where we clarified some central ideas regarding the span construction in higher categories.  Finally, the author thanks the University of California, Riverside where a first version of this paper was completed and where the author was supported as a graduate student in large part by the National Science Foundation under Grant No.\ 0653646.


\begin{thebibliography}{10}

\bibitem{aguilar} M.\ Aguiar and S.\ Mahajan, Monoidal functors,
species and Hopf algebras.  Available at
\href{http://www.math.tamu.edu/~maguiar/a.pdf}
{http://www.math.tamu.edu/$\sim$maguiar/a.pdf}

\bibitem{BaezCrans}
J.\ Baez and A.\ Crans, Higher-dimensional algebra VI: Lie 2-algebras, {\em Theory and Applications of Categories} {\bf 12} (2004), 492 �- 528. Also available as \href{http://arxiv.org/abs/math/0307263}{arXiv:math/0307263}.

\bibitem{BaezDolan:2001}
J.\ Baez and J.\ Dolan, From finite sets to Feynman diagrams, in {\em
Mathematics Unlimited---2001 and Beyond}, eds.\ Bj\"orn Engquist
and Wilfried Schmid, Springer, Berlin, 2001, pp.\ 29--50.  Also available as
\href{http://arxiv.org/abs/math/0004133}{arXiv:math/0004133}.

\bibitem{HDA7}
J.\ Baez, A.\ E.\ Hoffnung, and C.\ Walker, Higher Dimensional Algebra VII: Groupoidification, in {\em Theory and Applications of Categories}.  Available at\hfill\break
\href{http://arxiv.org/abs/0908.4305}{arXiv:0908.4305}.

\bibitem{HDA8}
J.\ Baez and A.\ E.\ Hoffnung, Higher Dimensional Algebra VIII: The Hecke Bicategory.  Available at \href{http://math.ucr.edu/home/baez/hecke.pdf}{http://math.ucr.edu/home/baez/hecke.pdf}.

\bibitem{BaezLangford}
J.\ Baez and L.\ Langford, Higher-dimensional algebra IV: 2-tangles, in {\em Adv. Math.} {\bf 180} (2003), 705 -� 764. Also available as \href{http://arxiv.org/abs/q-alg/9703033}
{arXiv:q-alg/9703033}.

\bibitem{BaezNeuchl}
J.\ Baez and M.\ Neuchl, Higher-dimensional algebra I: braided monoidal categories, in {\em Adv. Math.} {\bf 121} (1996), 196--244. Also available as \href{http://arxiv.org/abs/q-alg/9511013}
{arXiv:q-alg/9511013}.

\bibitem{Benabou}
J.\ B\'enabou, Introduction to Bicategories, in {\em Reports of the Midwest Category Seminar}, {\em Lecture Notes in Math.}, Springer-Verlag, Berlin, (1967), 1--77.

\bibitem{BFK}
J.\ Bernstein, I.\ Frenkel, and M.\ Khovanov, A categorification of the
Temperley-Lieb algebra and Schur quotients of $U(sl_2)$ via projective and
Zuckerman functors, \emph{Selecta Math.} {\bf 5} (1999), 199--241.

\bibitem{Brown}
K.\ Brown, {\sl Buildings}, Springer, Berlin, 1989.

\bibitem{Bump}
D.\ Bump, {\sl Lie groups}, Springer, New York, 2004.

\bibitem{Car}
S.\ M.\ Carmody, Cobordism Categories, PhD thesis, University of Cambridge, 1995.

\bibitem{CarterSaito}
J.\ S.\ Carter and M.\ Saito, {\sl Knotted Surfaces and Their Diagrams}, American Mathematical Society, Providence, 1998.

\bibitem{DayStreet}
B.\ Day and R.\ Street, Monoidal bicategories and Hopf algebroids, in {\em Adv. Math.} {\bf 129} (1997), 99�157.

\bibitem{ElKh}
B.\ Elias and M.\ Khovanov, Diagrammatics for Soergel Categories, 2009. Available as \href{http://arxiv.org/abs/0902.4700}{arXiv:0902.4700}.

\bibitem{FL} M.\ Fiore and T.\ Leinster, Objects of
categories as complex numbers, {\sl Adv.\ Math.\ }{\bf 190} (2005),
264--277.  Also available as
\href{http://arxiv.org/abs/math/0212377}{arXiv:math/0212377}.

\bibitem{For}
S.\ Forcey, Quotients of the multiplihedron as categorified associahedra, Available as \href{http://arxiv.org/PS_cache/arxiv/pdf/0803/0803.2694v4.pdf}{arXiv:0803.2694}.

\bibitem{FKS} I.\ B.\ Frenkel, M.\ Khovanov, and C.\ Stroppel, A categorification of finite-dimensional irreducible representations of quantum sl(2) and their tensor products, \emph{Selecta Math.} {\bf 12} (2006), 379--431.  Also available as \href{http://arxiv.org/abs/math/0511467}{arXiv:math/0511467}.

\bibitem{GPS}
R.\ Gordon, A.\ J.\ Powers, and R.\ Street, Coherence for Tricategories, {\em Mem. Amer. Math. Soc.} {\bf 117}, Providence, 1995.

\bibitem{Gurski}
N.\ Gurski, An algebraic theory of tricategories, PhD
thesis, University of Chicago, June 2006. Available as
\href{http://www.math.yale.edu/~mg622/tricats.pdf}{http://www.math.yale.edu/$\sim$mg622/tricats.pdf}

\bibitem{Hoffnung1}
A.\ E.\ Hoffnung, Spans in bicategories: A Semi-Strict Tetracategory, in progress.

\bibitem{Hoffnung2}
A.\ E.\ Hoffnung, On Enriched Bicategories, in progress.

\bibitem{Humphreys} J.\ Humphreys, {\sl Reflection Groups and
Coxeter Groups}, Cambridge U.\ Press, Cambridge, 1992.

\bibitem{Johnstone}
P.\ T.\ Johnstone, {\sl Sketches of an Elephant: A Topos Theory Compendium}, Oxford U.\ Press, Oxford, 2002.

\bibitem{Jones}
V.\ Jones, Hecke algebra representations of braid groups and link polynomials, {\sl Ann.\ Math.\ }{\bf 126} (1987), 335--388.

\bibitem{KaLu} D.\ Kazhdan and G.\ Lusztig, Representations of
Coxeter groups and Hecke algebras, \emph{Invent. Math.} {\bf 53}
(1979), 165--184.

\bibitem{KL}
M.\ Khovanov and A.\ Lauda, A diagrammatic approach to categorification of quantum groups I, {\em 	Represent. Theory} {\bf 13} (2009), 309--347.  Also available as \href{http://arxiv.org/abs/0803.4121}{arXiv:0803.4121}.

\bibitem{KV}
M.\ Kapranov and V.\ Voevodsky, 2-Categories and Zamolodchikov tetrahedra equations, in {\em Proc. Symp. Pure Math.} {\bf 56} Part 2 (1994), AMS, Providence, 177-�260.

\bibitem{KaLu} D.\ Kazhdan and G.\ Lusztig, Representations of
Coxeter groups and Hecke algebras, \emph{Invent. Math.} {\bf 53}, no.2 (1979), 165--184.

\bibitem{Kelly}
G.\ M.\ Kelly, {\sl Basic Concepts of Enriched Category Theory}, Cambridge University Press, Cambridge, 1982.

\bibitem{KePr}
T.\ Kenney and D.\ Pronk, Generalized span constructions, in progress.

\bibitem{Kh}
M.\ Khovanov, A categorification of the Jones polynomial, \textsl{Duke Math.\ J.\ }{\bf 101} (2000), 359--426.
Also available as \href{http://arxiv.org/abs/math/9908171}{arXiv:math/9908171}.

\bibitem{Kim}
M.\ Kim, A Lefschetz trace formula for equivariant cohomology, {\em Ann. Sci. \'Ecole Norm. Sup.} (4) {\bf 28} (1995), no. 6, 669 �- 688.

\bibitem{Leinster1}
T.\ Leinster, Basic Bicategories.  Available as \href{http://arxiv.org/abs/math/9810017}{arXiv:math/9810017}.

\bibitem{Leinster2}
T.\ Leinster, The Euler characteristic of a category, {\em Doc. Math.} {\bf 13} (2008), 21-�49. Also available as \href{http://arxiv.org/abs/math/0610260}{arXiv:math/0610260}.

\bibitem{MM}
S.\ Mac Lane and I.\ Moerdijk, {\sl Sheaves in Geometry and Logic: A First Introduction to Topos Theory}, Springer, Berlin, 1992.

\bibitem{MaStr}
V.\ Mazorchuk and C.\ Stroppel, Categorification of (induced) cell modules and the rough structure of generalised Verma modules, \emph{Advances in Mathematics} {\bf 219} (2008), 1363--1426.

\bibitem{Mc}
P.\ McCrudden, Balanced coalgebroids, in {\em Theory and Applications of Categories} {\bf 7}(6), 71-�147, 2000.

\bibitem{Moerdijk}
I.\ Moerdijk, Toposes and groupoids, in {\em Categorical algebra and its applications} (Louvain-La-Neuve, 1987),  280--298, Lecture Notes in Math., 1348, Springer, Berlin, 1988.

\bibitem{Rou} R.\ Rouquier, Categorification of the braid groups.  Available as \href{http://arxiv.org/abs/math/0409593}{arXiv:math/0409593}.

\bibitem{Soe}  W.\ Soergel, The combinatorics of Harish-Chandra bimodules, \emph{Journal Reine Angew. Math.} {\bf 429} (1992), 49--74.

\bibitem{Str1} C.\ Stroppel, Category $\mathcal{O}$: Gradings and translation functors, \emph{J. Algebra} {\bf 268} (2003), 301--326.

\bibitem{Str2} C.\ Stroppel, Categorification of the Temperley�Lieb category, tangles, and cobordisms via projective functors, \emph{Duke
Math. J.} {\bf 126} (3) (2005) 547--596.

\bibitem{Weber} M.\ Weber, Strict $2$-toposes.  Available as \href{http://arxiv.org/abs/math/0606393}{arXiv:math/0606393}.

\bibitem{WW} B.\ Webster and G.\ Williamson, A geometric model for Hochschild homology of Soergel bimodules, \emph{Geometry and Topology} {\bf 12} (2008), 1243--1263.  Also available as \href{http://arxiv.org/abs/0707.2003}{arXiv:0707.2003}.

\bibitem{Weinstein}
A.\ Weinstein, The Volume of a Differentiable Stack, {\em proceedings of Poisson 2008 (Lausanne, July 2008)}.  Also available as \href{http://arxiv.org/abs/0809.2130}{arXiv:0809.2130}.
\end{thebibliography}
\end{document}